\documentclass[preprint]{elsarticle}

\usepackage{amsmath,amssymb,mathdots,mathrsfs}
\usepackage{amsthm}
\usepackage{comment}
\usepackage{pgfplotstable}
\usepackage{times}

\usepackage[colorlinks,bookmarksopen,bookmarksnumbered,citecolor=red,urlcolor=red, unicode]{hyperref}
\usepackage{url}

\def\scrH{\mathscr{H}}

\def\scrS{\mathscr{S}}

\def\calS{\mathcal{S}}

\def\calX{\mathcal{X}}

\def\bfc{\mathbf{c}}

\def\bfg{\mathbf{g}}

\def\bfk{\mathbf{k}}
\def\bfl{\mathbf{l}}

\def\bfn{\mathbf{n}}

\def\bfp{\mathbf{p}}
\def\bfq{\mathbf{q}}

\def\bfu{\mathbf{u}}

\def\bfw{\mathbf{w}}

\def\bfone{\mathbf{1}} 
\def\bfzero{\mathbf{0}} 

\def\bbC{\mathbb{C}}

\def\bbF{\mathbb{F}}

\def\bbN{\mathbb{N}}

\def\bbR{\mathbb{R}}

\def\bbT{\mathbb{T}}

\def\rmF{\mathrm{F}}

\def\rmF{\mathrm{F}}

\newcommand{\defeq}{\mathrel{:=}}
\newcommand{\mvec}{\mathop{\text{vec}}}
\newcommand{\rrank}{\mathop{\text{rank}}}

\newcommand{\minim}{\mathop{\text{\;minimize\;}}}
\newcommand{\sto}{\mathop{\text{\;subject\ to\;}}}
\newcommand{\diag}{\mathop{\mathrm{diag}}} 
\newcommand{\vcol}{\mathop{\mathrm{col}}} 
\def\blkdiag{\mathop{\mathrm{blkdiag}}} 

\newcommand{\bmx}{\begin{bmatrix}}
\newcommand{\emx}{\end{bmatrix}}
\newcommand{\bsm}{\left[\begin{smallmatrix}}
\newcommand{\esm}{\end{smallmatrix}\right]}

\newcommand{\ie}{\emph{i.e.}}

\def\rowvec#1#2{[\,#1 \; \cdots \; #2\,]}
\def\pspace#1{\mathcal{P}_{#1}} 

\def\multmat#1#2{\mathbf{M}_{#2}({#1})}

\def\pco#1#2{#1_{#2}}
\def\pno#1#2{#1^{(#2)}}
\def\phatno#1#2{{\widehat{#1}}^{(#2)}}

\def\sylv#1{\mathcal{S}_{#1}} 

\def\elsum#1{\Sigma(#1)}

\def\gcddist{\mathop{\mathrm{dist}}}
\def\pspacen{\pspace{\bfn}}
\def\pspacecd#1{\mathcal{F}_{#1}}
\def\pspacegcd#1{\mathcal{G}_{#1}}
\def\matmultmat#1#2{\mathcal{A}_{#1,#2}}

\def\lsfun#1{g^{(#1)}} 
\def\optfun#1{f^{(#1)}_*} 
\def\lsfunopt#1{g^{(#1)}_*}

\newtheorem{theorem}{Theorem}[section]
\newtheorem{algorithm}{Algorithm}[section] 
\newtheorem{prob}{Problem}[section] 
\newtheorem{example}{Example}[section]
\newtheorem{corollary}[theorem]{Corollary}
\newtheorem{lem}[theorem]{Lemma} 
\newtheorem{prop}[theorem]{Proposition} 
 
\newtheorem{rem}[theorem]{Remark} 
\newtheorem{definition}[theorem]{Definition}

\newenvironment{pf}%
	       {{\bf Proof}.\ }{\hfill$\Box$}

\begin{document}
\author[g1,g2]{Konstantin Usevich\corref{cor}}\ead{konstantin.usevich@gipsa-lab.grenoble-inp.fr}
\author[b]{Ivan Markovsky}\ead{ivan.markovsky@vub.ac.be}
\cortext[cor]{Corresponding author}
\address[g1]{Univ. Grenoble Alpes, GIPSA-Lab, F-38000 Grenoble, France}
\address[g2]{CNRS, GIPSA-Lab, F-38000 Grenoble, France}
\address[b]{Department ELEC, Vrije Universiteit Brussel (VUB), Pleinlaan 2,  B-1050 Brussels, Belgium}
\title{Variable projection methods for approximate (greatest) common divisor computations}
\date{}

\begin{abstract}

We consider the problem of finding for a given $N$-tuple of polynomials (real or complex) the closest $N$-tuple that has a common divisor of degree at least $d$. Extended weighted Euclidean seminorm of the coefficients is used as a measure of closeness. Two equivalent representations of the problem are considered: (i) direct parameterization over the common divisors and quotients (image representation), and (ii) Sylvester low-rank approximation (kernel representation). We use the duality between least-squares and least-norm problems to show that (i) and (ii) are closely related to mosaic Hankel low-rank approximation. This allows us to apply to the approximate common divisor problem recent results on complexity and accuracy of computations for mosaic Hankel low-rank approximation. We develop optimization methods based on the variable projection principle both for image and kernel representation. These methods have linear  complexity in the degrees of the polynomials for small and large $d$. We provide a software implementation of the developed methods, which is based on a software package for structured low-rank approximation. 
\end{abstract}

\begin{keyword}
approximate GCD; structured low-rank approximation; variable projection; mosaic Hankel matrices; least squares  problem; weighted 2-norm
\end{keyword}


\maketitle

\section{Introduction}
The problem of computing a greatest common divisor (GCD) of polynomials with real or complex coefficients appears in many applications: signal processing and system identification \cite{Agrawal.etal04SP-Common,Gaubitch.etal05conf-Adaptive}, computer-aided geometric design \cite{Emiris.etal13CD-Sparse}, blind image deblurring \cite{Li.etal10conf-Blind}, control of linear systems \cite{Khare.etal10conf-Real,uncontr} and approximate factorization of polynomials \cite{Zeng08-ApaTools}.
But, as noted in \cite{Pan01IaC-Computation}, ``computation of polynomial GCDs is an excellent example of numerically ill-posed problems''.
Indeed, for any set of polynomials with non-trivial GCD, a generic perturbation of the coefficients makes them coprime.
In the aforementioned applications, perturbations appear naturally due to limited precision of the numerical representation of the coefficients, or due to measurement errors.
These reasons make it inevitable to use the notion of approximate greatest common divisor (AGCD).

There is a vast literature on the topic of AGCD (see the list of references), starting with different definitions of AGCD and finishing with different computational methods.
Nevertheless, two main optimization-based formulations of AGCD are predominant.
In what follows, we give these formulations for the case of two polynomials.

The first commonly accepted problem formulation is the problem of finding the so-called $\varepsilon$-GCD (see \cite[Def. 2.2]{Rupprecht99JPAA-algorithm}, \cite[Def 1.1]{Bini.Boito10-fast}, \cite[Eqn. (1)-(2)]{Pan01IaC-Computation}).
\begin{prob}[$\varepsilon$-GCD]\label{prob:agcd}
Given polynomials $p(z), q(z)$, and a threshold $\varepsilon$, find polynomials $\widehat{p}^{*}(z), \widehat{q}^{*}(z)$ with a GCD $\widehat{h}^{*}(z) = \gcd(\widehat{p}^{*}(z), \widehat{q}^{*}(z))$ that are solutions to 
\begin{equation}\label{eq:problem_agcd}\tag{AGCD}
\max_{\widehat{p}(z), \widehat{q}(z)} \deg \gcd(\widehat{p}(z), \widehat{q}(z)) \;\sto\; 
\gcddist\big((p,q),(\widehat{p}, \widehat{q})\big) \le \varepsilon,
\end{equation}
where $\gcddist(\cdot,\cdot)$ is some distance measure for the pairs of polynomials. 
\end{prob}
The polynomial $\widehat{h}^{*}(z)$ is conventionally called an $\varepsilon$-GCD. The distance in the definition of the $\varepsilon$-GCD is typically of the form 
\[\gcddist \big((p,q),(\widehat{p}, \widehat{q})\big) = \|(p-\widehat{p},q-\widehat{q})\|,\] 
where $\|\cdot\|$ is some norm.
Various norms are used in the literature: the $\ell_2$-norm (\cite{Cheze.etal11TCS-subdivision,Agrawal.etal04SP-Common,Corless.etal95conf-singular,Kaltofen.etal07-Structured}), the  $\ell_\infty$-norm \cite{Hitz98PhD-Efficient} and mixed $\ell_2$/$\ell_\infty$-norm in \cite[Def. 2.2]{Rupprecht99JPAA-algorithm} and \cite[Def 1.1]{Bini.Boito10-fast} (i.e., $\|(p,q)\| = \max(\|p\|_2, \|q\|_2)$).

The second core problem formulation \cite[Prob. 1.1]{Kaltofen.etal07-Structured} is a dual problem to \eqref{eq:problem_agcd}.
\begin{prob}\label{prob:acd}
Given $p(z), q(z)$ and a number $d$, find $\widehat{p}(z), \widehat{q}(z)$  that are solutions to 
\begin{equation}\label{eq:problem_acd}\tag{ACD}
\min_{\widehat{p}(z), \widehat{q}(z)} \gcddist \big((p,q),(\widehat{p}, \widehat{q})\big) 
\sto \deg \gcd(\widehat{p}(z), \widehat{q}(z)) \ge d.
\end{equation}
\end{prob}
Problem~\ref{prob:acd} appears in many contexts.
First, it is often used in a combination with Problem~\ref{prob:agcd}: since Problem~\ref{prob:agcd} typically has infinite number of optimal solutions $(\widehat{p}^{*},\widehat{q}^{*})$, the closest pair of polynomials to the given ones is of interest.
Thus Problem~\ref{prob:acd} is often referred to as \emph{refinement} in the AGCD literature \cite[\S 3.2]{Bini.Boito10-fast}, \cite{Zeng11-Numerical}.

\begin{figure}[ht]%
\centering%
\includegraphics[width=0.7\textwidth]{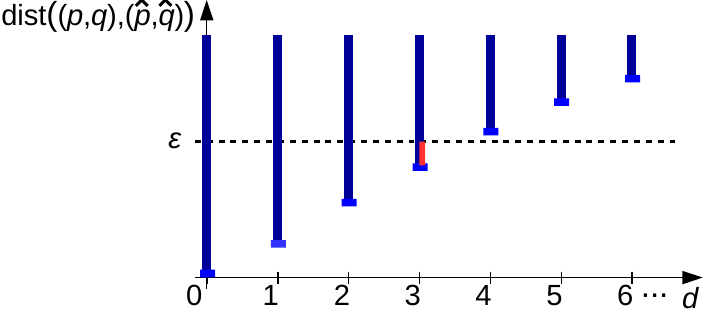}%
\caption{Dark blue: feasible set for \eqref{eq:problem_agcd} and \eqref{eq:problem_acd}; light blue: optimal points for \eqref{eq:problem_acd}; red:  optimal points for \eqref{eq:problem_agcd} for a given error $\varepsilon$ and degree $d$. }%
\label{fig:pareto}
\end{figure}

Second, as noted in \cite{Rupprecht99JPAA-algorithm}, being able to solve Problem~\ref{prob:acd},  gives a solution to Problem~\ref{prob:agcd}. Indeed, if the minimum value of  \eqref{eq:problem_acd} for $d = d^*$ is less than or equal to a given $\varepsilon$ and the minimum value of  \eqref{eq:problem_acd} for $d = d^*+1$ is strictly greater than $\varepsilon$, then $d^*$ is the solution of \eqref{eq:problem_agcd}.
(For more details, see the discussion after \cite[Def. 2.3]{Rupprecht99JPAA-algorithm}.)
This fact is illustrated in Fig.~\ref{fig:pareto}, where the feasible set  for \eqref{eq:problem_agcd} and \eqref{eq:problem_acd} is shown (i.e., the values of $d$ and $\gcddist \big((p,q),(\widehat{p}, \widehat{q})\big)$, for which exists a pair $(\widehat{p}, \widehat{q})$ that satisfies  $\deg\gcd(\widehat{p}, \widehat{q}) \ge d$). 
The solutions of \eqref{eq:problem_acd} correspond to the lowest points in each vertical line.
The optimal solutions of \eqref{eq:problem_agcd} correspond to red segment in Fig.~\ref{fig:pareto} (the points on the rightmost vertical line that intersect the threshold horizontal line).
Thus, Problem~\ref{prob:agcd} can be solved by solving Problem~\ref{prob:acd} for all possible $d$ (or by using, for example, bisection over $d$).

Finally, Problem~\ref{prob:acd} appears when we know an \emph{a priori} bound on the degree of the GCD is given, which is a reasonable assumption in some applications \cite{Agrawal.etal04SP-Common}, \cite{Markovsky.VanHuffel06conf-algorithm}.
Another common example is  the problem of finding the nearest non-coprime polynomials \cite{Cheze.etal11TCS-subdivision, uncontr}, which corresponds to $d=1$.
In this paper, we focus on Problem~\ref{prob:acd}.

\subsection{Previous works}
Two main approaches to Problem~\ref{prob:acd} can be identified in the literature, namely the  direct parameterization approach (referred to as \emph{image representation} in this paper) and the structured low-rank approximation (SLRA) approach (also referred to as \emph{kernel representation} in this paper).
Most of the algorithms were proposed for minimizing the weighted Euclidean distance.

The image representation approach is based on the direct representation of the polynomials as a product of a common factor and quotient polynomials, i.e., the cost function $f(\widehat{h}, \widehat{u}, \widehat{v}) \defeq  \gcddist((p,q), (\widehat{h}\widehat{u},\widehat{h}\widehat{v}))$ is minimized over all candidate common divisors $\widehat{h}$ (of degree $d$) and candidate quotient polynomials $\widehat{u}$, $\widehat{v}$ (of degree $n-d$).
The image representation approach is used as early as \cite{Corless.etal95conf-singular}.
However, the size of the search space makes minimization of $f$ expensive for general-purpose optimization routines. 
One of the main ways to reduce the complexity is elimination  of variables, which is also named  variable projection in the context of nonlinear least squares problems \cite{GolubPereyra73SJoNA-Differentiation}. 

In short, the variable projection principle is: for a fixed $\widehat{h}$, the minimization of $f$ with respect to the other parameters is a linear least squares problem; thus other parameters can be eliminated, and a function with smaller number of parameters $f(\widehat{h})$ can be minimized.
In the special case of $d=1$, as shown in \cite{Karmarkar.Lakshman98JoSC-Approximate}, the minimum of $f(\widehat{h})$ can be computed by minimizing a univariate polynomial (for real $p$ and $q$) or bivariate real polynomial (for complex $p$ and $q$).
For the latter case a certified algorithm was presented in \cite{Cheze.etal11TCS-subdivision}.
If $d > 1$, but $d$ is small, as shown in \cite{Corless.etal95conf-singular,Pan01IaC-Computation}, $f(\widehat{h})$ can be computed efficiently (in linear time in the degree of the polynomials, if $d$ is small). 
Later (but independently) in \cite{Markovsky.VanHuffel06conf-algorithm}, it was shown that the first derivative of $f(\widehat{h})$ can be evaluated with the same complexity.
Independently, in \cite{Stoica.Soderstrom97A-Common,Agrawal.etal04SP-Common}, elimination of $\widehat{h}$ (instead of $\widehat{u}$ and $\widehat{v}$) was proposed.
Finally, the variable projection was also implicitly used  in \cite{Zhi.etal04JJoIaAM-Hybrid,Li.Zhi13TCS-Computing} for symbolic computation of nearest singular polynomial and semi-definite programming relaxations of  the AGCD problem \cite{Li.etal08TCS-Approximate,Kaltofen.etal08conf-Exact}.
There are other developments and extensions for the image representation.
For example, in \cite[\S 3.2]{Bini.Boito10-fast} it was shown that the Gauss-Newton step can be performed in quadratic time in the degrees of the polynomials.

Another popular approach to Problem~\ref{prob:acd} is the SLRA approach (kernel representation approach), which consists in reformulating Problem~\ref{prob:acd} as a problem of approximating a given structured matrix by a structured matrix of low rank (an SLRA problem \cite{Markovsky08A-Structured}).
This reformulation is possible since the constraint on the GCD degree can be rewritten as a rank constraint on a structured matrix.
For two polynomials, this is a Sylvester or a Sylvester subresultant matrix \cite{Zeng11-Numerical}; for several polynomials there are various generalizations of the Sylvester  structure \cite{Kaltofen.etal06conf-Approximate,Karcanias.etal06CMwA-Approximate} (see Section \ref{sec:slra} for more details).

Concerning algorithms for SLRA, the following methods were used in the context of the AGCD problem:
 {\em structured total least norm (STLN)\/} and its improvements \cite{Kaltofen.etal07-Structured,Li.etal05JJSSaAC-Fast,Winkler.Allan08JoCaAM-Structured,Winkler.Hasan13JoCaAM-improved}, 
{\em Riemannian SVD\/} \cite{Botting.etal05conf-Using}, 
 {\em gradient projection\/} \cite{Terui13TCS-GPGCD}, 
 {\em alternating least squares with penalization\/} \cite{Ishteva.etal13-Regularized}, 
 {\em variable projection\/} \cite{Markovsky.Usevich13JCAM-Software}, 
{\em Newton-like alternating projection algorithms\/} \cite{Schost.Spaenlehauer13-Quadratically}.
In \cite{Ottaviani.etal13-Exact}, a step toward global optimization was made by the authors who proposed to compute the number of complex critical points for the optimization problem, using symbolic computations.

\subsection{Contribution and structure of this paper}

In this paper, we consider generalization of Problem~\ref{prob:acd} to many polynomials.
The main contributions of this paper are connections between image/kernel representation approaches to Problem~\ref{prob:acd} and structured low-rank approximation of mosaic-Hankel matrices.
First, we show that the generalized Sylvester subresultant low-rank matrix approximation can be reduced to mosaic Hankel SLRA.
Second, we show that the cost function in variable projection methods for AGCD in image representation has the same structure as the cost function in the variable projection method for mosaic Hankel SLRA \cite{Usevich.Markovsky13JCAM-Variable}.
These connections allow us to use, with small modification, the efficient algorithms developed in \cite{Usevich.Markovsky13JCAM-Variable}. 
The algorithms have proven computational complexity and can handle real and complex polynomials.
As a side result, we show that minimizing the relative distance between tuples of polynomials is equivalent to minimizing a distance based on angles between polynomials.

This structure of the paper is as follows. Sections~\ref{sec:prel}--\ref{sec:mosaic_slra} contain known results or their minor improvements. Sections~\ref{sec:results}--\ref{sec:exper} contain the main results and experiments. Section~\ref{sec:prel} contains the necessary background and a formal statement of an analogue of Problem~\ref{prob:acd} for many polynomials;
we also introduce the spaces of homogeneous polynomials (polynomials with possible infinite roots) and operations with them, which are key ingredients of this paper.
In Section~\ref{sec:nls} we review the image representation approach and the variable projection principle.
In Section~\ref{sec:slra}, we review the structured low-rank approximation  (kernel representation) approaches adapted to our problem statement.
In Section~\ref{sec:mosaic_slra}, we recall the mosaic Hankel structure and the results on variable projection methods of the corresponding SLRA problem.
In Section~\ref{sec:results}, we present the main results of the paper.
In Section~\ref{sec:exper}, we provide numerical experiments that include comparison with the state-of-the-art methods.
The methods developed in this paper are implemented in MATLAB and  are based on the SLRA package \cite{SLRA} described in \cite{Markovsky.Usevich13JCAM-Software}. The source code of the methods and experiments is publicly available at \url{http://github.com/slra/slra}. 

\section{Main notation and the approximate common divisor problem}\label{sec:prel}

\subsection{Polynomials with possible roots at infinity}
Let $\bbF[z]$ denote the set of univariate polynomials over the field $\bbF$ (where $\bbF$ is $\bbC$ or $\bbR$).
Let $\pspace{n} \subset \bbF[z]$ be the set of  polynomials of degree at most $n \ge 0$, i.e.
\begin{equation}
\pspace{n} := \{ p_0 + p_1 z + \cdots + p_n z^n  \; | \; p_j \in \bbF \} \subset \bbF[z].
\label{eq:def_pspace}
\end{equation}
Then the space $\pspace{n}$  is isomorphic to $\bbF^{n+1}$ through the correspondence
\begin{equation}\label{eq:poly_vec_isomorphism}
p(z) = p_0 + p_1 z + \cdots + p_n z^n \in \pspace{n} \leftrightarrow p = \rowvec{p_0}{p_n}^{\top} \in \bbF^{n+1}.
\end{equation}
With some possible abuse of notation, we will use shorthand notation $p \in \pspace{n}$.

\begin{rem}
The leading coefficient $p_n$ may be equal to $0$. In this case, we will say that the polynomial $p(z)$ has the  root $\infty$. 
By multiplicity of the root $\infty$ we denote the maximal number of consecutive zero leading coefficients.
Thus every polynomial in $\pspace{n}\setminus \{0\}$ has exactly $n$ roots  in $\bbC \cup \{\infty\}$ (the Riemannian sphere).
\end{rem}

\begin{example}
The following polynomial has two simple roots ($1$ and $2$) and a double root $\infty$:
\[
0\cdot z^4 + 0\cdot z^3 + z^2 - 3z +2 \in \pspace{4}
\] 

\end{example}

\begin{rem}
In this paper, we call the elements of $\pspace{n}$ \textit{homogeneous polynomials}, since  $\pspace{n}$ can be viewed as the space of bivariate homogeneous polynomials.
\end{rem}

\subsection{Multiplication and division of homogeneous polynomials}\label{sec:pol_mult}
The  multiplication of polynomials is defined as $(\pno{p}{1} \cdot \pno{p}{2}) (z) = \pno{p}{1}(z) \pno{p}{2} (z)$ (acting as $\pspace{n_1} \times \pspace{n_2} \to \pspace{n_1+n_2}$). It has the following matrix representation:
\[
\pno{p}{1} \cdot \pno{p}{2} = \multmat{\pno{p}{1}}{n_2} \pno{p}{2} = \multmat{\pno{p}{2}}{n_1} \pno{p}{1},
\]
where $\pno{p}{1} \in \bbF^{n_1+1}$,$\pno{p}{2} \in \bbF^{n_2+1}$, and $\multmat{h}{m}$ is the \textit{multiplication matrix}  by  $h \in \pspace{d}$:
\[
\multmat{h}{m} \defeq
\bmx 
\pco{h}{0} &        &            \\ 
\vdots     & \ddots &            \\
\pco{h}{d} &        & \pco{h}{0} \\
           & \ddots & \vdots     \\  
           &        & \pco{h}{d}    
\emx \in \bbF^{(m+d+1) \times (m+1)},
\]
which is  a rectangular Toeplitz matrix,
where the blank triangular parts stand for zeros. 

For $0 \le d \le n$, we say that a polynomial $h \in \pspace{d} \setminus \{\bfzero\}$ \textit{divides} a polynomial $p \in \pspace{n}$ (or $h$ is a \textit{divisor} of $p$), if there exists a polynomial $q \in \pspace{n-d}$ such that $p = q\cdot h$. In particular, this definition includes the following special cases.
\begin{itemize}
\item All  $h \in \pspace{d} \setminus \{\bfzero\}$, $0 \le d \le n$, are the divisors of the zero polynomial $\bfzero \in \pspace{n}$.
\item A nonzero polynomial of zero degree $h \in \pspace{0} \setminus \{ \bfzero \}$ is a divisor of any polynomial.
\end{itemize}

The notion of divisor in the spaces $\pspace{n}$ differs from the notion of divisors for ordinary polynomials, due to possible presence of the $\infty$ roots.

\begin{example}
Consider two polynomials
\[
p(z) = z^4 - 3z^3 +2z^2 \in \pspace{4} \mbox{\; and \;} h(z) = 0\cdot z^2 + z-1 \in \pspace{2}.
\]
Although $(z-1)$ divides of $p$, the polynomial $h$ does not, because it has the root $\infty$.
\end{example}

\subsection{N-tuples of polynomials and common divisors}\label{sec:n_tuples}
Let $\bfn = \rowvec{n_1}{n_N}^{\top} \in \bbN^{N}$ be a vector of fixed degrees, $n_{min} \defeq \min n_k$, $n \defeq \sum_{k=1}^N n_k$, and denote by $\pspacen \defeq \pspace{n_1} \times  \cdots \times \pspace{n_N}$ the set of $N$-tuples of polynomials with these degrees. We also adopt the notation $\bfp = (\pno{p}{1},\ldots,\pno{p}{N})$, $\pno{p}{k} \in \pspace{n_k}$, for the elements of the $N$-tuples.
With some possible ambiguity of notation (as in 
\eqref{eq:poly_vec_isomorphism}), we use the same letter for the $N$-tuple $\bfp \in \pspace{\bfn}$, and for the stacked vector
\[
\bfp = \vcol(\pno{p}{1},\ldots,\pno{p}{N}) \in \bbF^{n_1+\cdots+n_N+N}.
\]

Next, we introduce the operation of multiplication of an $N$-tuple $\bfg = (\pno{g}{1}, \ldots, \pno{g}{N}) \in \pspace{\bfn-d} = \pspace{n_1-d} \times  \cdots \times \pspace{n_N-d}$ by a polynomial $h \in \pspace{d}$ as follows:
\[
\bfg \cdot h \defeq (\pno{g}{1} \cdot h, \ldots, \pno{g}{N} \cdot h).
\]

\begin{definition}
The polynomials $\bfp = (\pno{p}{1}, \ldots, \pno{p}{N}) \in \pspacen$ have a \textit{common divisor}  $h \in \pspace{d} \setminus \{\bfzero\}$, if $h$ divides all the polynomials $\pno{p}{k}$. The polynomial $h$ is called a \emph{greatest common divisor (GCD)} if there are no common divisors in $\pspace{d'}$ for any $d' > d$.
\end{definition}

 Since  $1 \in \pspace{0}$ is a divisor of any polynomial, a GCD always exists (but it is not unique). We denote by $\deg \gcd \bfp$ the degree of GCDs, and denote by $\gcd \bfp$ the set of all GCD (which is a subset of $\pspace{\deg \gcd \bfp}$). In particular, if all $\bfp = 0$, then  
\[
\gcd \bfp = \pspace{n_{min}}\setminus \{\bfzero\}.
\]
Otherwise, $\gcd \bfp$ is a punctured one-dimensional linear subspace of $\pspace{\deg \gcd \bfp}$
\[
\gcd \bfp = \{\alpha h : \alpha \in \bbF \setminus \{0\},\quad h \; \mbox{is a GCD of } \pno{p}{1}, \ldots, \pno{p}{N} \}.
\]

\subsection{The approximate common divisor problem statement}\label{sec:prob_stat}
Define the set of $N$-tuples that have a GCD of degree at least $d \ge 0$ as follows
\begin{equation}
\pspacegcd{d} \defeq 
\left\{\widehat{\bfp} = (\pno{\widehat{p}}{1}, \ldots, \pno{\widehat{p}}{N}) \in \pspace{\bfn} \;|\;
\deg\gcd (\widehat{\bfp}) \ge d \right\}.\label{eq:mod_class2}
\end{equation}
Finally, assume that $\pspacen$ is equipped with a distance $\gcddist(\cdot, \cdot)$, which is continuous in the Euclidean topology.
We formulate the generalization of Problem~\ref{prob:acd} as follows.
\begin{prob}[Approximate GCD with bounded degree]\label{pro:acd2}
Given  $\bfp = (\pno{p}{1} , \ldots, \pno{p}{N}) \in \pspacen$ and  $d:0 \le d \le n_{min}$, find the distance
\begin{equation}
\gcddist\left( \bfp, \pspacegcd{d} \right) \defeq
\min_{\widehat{\bfp} \in \pspacegcd{d}}  \quad
\gcddist\left(\bfp, \widehat{\bfp}\right).\label{eq:acd2}
\end{equation}
\end{prob}
Note that 
Problem~\ref{pro:acd2} well-posed, since the  set $\pspacegcd{d}$ is closed in the Euclidean topology, by Lemma~\ref{lem:mod_class_relation} (see \ref{sec:acd_basic}).

\begin{rem}
The set $\pspacegcd{d}$ is closed, in particular, thanks to the fact that we use homogeneous polynomials and to the special definition of the GCD in Section~\ref{sec:n_tuples}.
\end{rem}

\subsection{Weigthed norm, missing and fixed values}\label{sec:wnorm}
In this paper,  we use the following distance:
\begin{equation}\label{eq:wdist}
\gcddist(\bfp,\bfq) \defeq \|\bfp - \bfq\|^2_{\bfw} = \sum\limits_{k=1}^N \|\pno{p}{k} - \pno{q}{k}\|_{\pno{w}{k}}^2,
\end{equation}
where $\bfw = (\pno{w}{1},\ldots,\pno{w}{N})$ is a tuple of weight vectors $\pno{w}{k} \in [0,+\infty]^{n_k+1}$,  and
\begin{equation}\label{eq:wnorm}
 \|p\|^2_w \defeq \sum_{j=0}^{n} w_j p_j \overline{p_j},
\end{equation}
is the weighted extended semi-norm on $\pspacen$ (see also \cite{Markovsky.Usevich13JCAM-Software}).
If $w_j \in (0,\infty)$, then $\|\cdot\|^2_{\bfw}$ is a standard weighted $\ell_2$-norm. The $0$ and $\infty$ weights have a special meaning:
\begin{itemize}
\item We formally assume $0 \cdot \infty = 0$ in \eqref{eq:wnorm} and require the distance to be finite in \eqref{eq:acd2}. Hence, a weight $\pno{w}{k}_j = \infty$  is equivalent to an equality constraint $\pno{p}{k}_j = \pno{\widehat{p}}{k}_j$. For example,  monicity constraints  can be imposed using  $\infty$ weights. 
\item If a  weight $\pno{w}{k}_j = 0$ is present,  the solution $\widehat{\bfp}$ does not depend on $\pno{p}{k}_j$.  Hence, we may assume that the coefficient is undefined (the case of missing data  \cite{Markovsky.Usevich13SJMAA-Structured}).
\end{itemize}

\section{Image representation and variable projection}\label{sec:nls}
\subsection{Image representation}
In this approach, the set $\pspacegcd{d}$ is replaced by the set (of candidate factorizations)
\begin{equation}\label{eq:mod_class}
\pspacecd{d} \defeq 
\left\{(\pno{g}{1} h, \ldots, \pno{g}{N} h) \;|\; (\pno{g}{1}, \ldots, \pno{g}{N}) \in \pspace{\bfn-d}, h \in \pspace{d} \setminus \{\bfzero\} \right\},
\end{equation}
where $h$ is a candidate common divisor, and $\pno{g}{k}$ are candidate quotient polynomials. 

We refer to \eqref{eq:mod_class} as \emph{image representation}, since  $\pspacecd{d}$ is the image of the map
\[
\begin{split}
  \pspace{\bfn-d} \times (\pspace{d} \setminus \{ \bfzero \}) & \to   \pspace{\bfn},  \\
 (\widehat{\bfg},\widehat{h}) & \mapsto \widehat{\bfg} \cdot \widehat{h}.
\end{split}
\]
\begin{rem}\label{rem:two_classes}
For complex polynomials, the sets $\pspacegcd{d}$ and $\pspacecd{d}$ coincide.
But for real polynomials we have $\pspacecd{d} \subseteq \pspacegcd{d}$, and the equality is not always satisfied.
The precise relation between  $\pspacegcd{d}$ and $\pspacecd{d}$ is given in Lemma~\ref{lem:mod_class_relation} in \ref{sec:acd_basic}.
\end{rem}
Then an analogue of Problem~\ref{pro:acd2} is  formulated as follows.
\begin{prob}\label{pro:agcd_nls}
Given $\bfp$, find $\widehat{h}$ and $\widehat{\bfg} = (\pno{g}{1}, \ldots, \pno{g}{N})$ that are solutions to
\begin{equation}
\minim_{\footnotesize\begin{array}{c}\widehat{\bfg} \in \pspace{\bfn-d}, 
\widehat{h} \in \pspace{d}\end{array}}  
\| \widehat{\bfg} \cdot \widehat{h}  - \bfp\|_{\bfw}^{2}.
\label{eq:agcd_nls}
\end{equation}
\end{prob}
\begin{rem}
We include the zero polynomial $\widehat{h}$ in the search space of \eqref{eq:agcd_nls} (compared to the definition in \eqref{eq:mod_class}), since the zero tuple $\bfzero \in \pspace{\bfn}$ belongs to $\pspacecd{d}$ by definition.
\end{rem}
\begin{rem}
By Remark~\ref{rem:two_classes}, Problems~\ref{pro:acd2} and \ref{pro:agcd_nls} are not equivalent.
\end{rem}

\subsection{Variable projection methods}\label{sec:varpro_methods}
Denote the cost function in \eqref{eq:agcd_nls} as
\begin{equation}
{f}(\widehat{h},\widehat{\bfg}) = \| \widehat{\bfg} \cdot \widehat{h}  - \bfp\|_{\bfw}^{2}.
\label{eq:cost_sec_varpro}
\end{equation}
The cost function \eqref{eq:cost_sec_varpro} is a  nonlinear least squares problem, for which the variable projection principle \cite{GolubPereyra73SJoNA-Differentiation} (based on elimination of variables) can be applied.

The variable projection principle \cite{GolubPereyra73SJoNA-Differentiation} is based on the fact that for one fixed  variable (either  $\widehat{h}$ or $\widehat{\bfg}$), minimization of \eqref{eq:cost_sec_varpro} is a linear least squares problem and has a closed form solution. This principle is further explained on each of the examples.
\begin{example}[Variable projection with respect to a common divisor]\label{ex:varpro_h}
In this case, the problem \eqref{pro:agcd_nls} is rewritten as the following double minimization problem
\begin{align}
&\minim_{\widehat{h} \in \pspace{d} \setminus \{0\}} f_1(\widehat{h}), \quad \mbox{where}  \label{eq:outer} \\
&f_1(\widehat{h}) \defeq \min_{\widehat{\bfg}\in \pspace{\bfn-d}} f(\widehat{h},\widehat{\bfg}).\nonumber
\end{align}
We note that computing $f_1(\widehat{h})$ can be rewritten as a linear least squares problem, 
\[
f_1(\widehat{h})  = \min_{\widehat{\bfg}\in \pspace{\bfn-d}}
 \sum\limits_{k=1}^N 
 \|\multmat{\widehat{h}}{n_k-d} \pno{\widehat{g}}{k}- \pno{p}{k}\|^2_{\pno{w}{k}},
\]
which has a closed form solution (see \ref{sec:ls_ln}).
Thus  in \eqref{eq:outer}  $\widehat{\bfg}$ is eliminated.
\end{example}

\begin{example}[Variable projection with respect to quotient polynomials]\label{ex:varpro_bfg}
In this case, the problem \eqref{pro:agcd_nls} is rewritten as the following double minimization problem
\begin{align}
&\minim_{\widehat{\bfg} \in \pspace{\bfn-d} \setminus \{0\}} f_2(\widehat{\bfg}), \quad \mbox{where}  \label{eq:outer_cofe} \\
&f_2(\widehat{\bfg}) \defeq \min_{\widehat{h}\in \pspace{d}}f(\widehat{h},\widehat{\bfg}). 
\nonumber
\end{align}
As in Example~\ref{ex:varpro_h}, the function $f_2$ has the form
\begin{equation}\label{eq:inner_cofe} 
f_2(\widehat{\bfg}) = \min_{\widehat{h} \in \pspace{d}} \sum\limits_{k=1}^N 
 \|\multmat{\pno{\widehat{g}}{k}}{d} \widehat{h} - \pno{p}{k}\|^2_{\pno{w}{k}},
\end{equation}
thus the variable $\widehat{h}$  is eliminated from the problem \eqref{eq:outer_cofe}.
\end{example}

\begin{rem}[On algorithms]\label{rem:algorithms}
After eliminated of variables, for the reduced cost function ($f: \calX \to \bbR$)  we can apply conventional smooth optimization methods, such as:
\begin{itemize}
\item gradient-based methods, which require evaluation of the cost function $f(x)$ and its gradient $\nabla f(x)$ (including quasi-Newton methods, for example, BFGS \cite{Nocedal.Wright06-Numerical});
\item if the cost function can be represented as a sum of squares, i.e.,
\[
f(x) = \|g(x)\|^2_2,
\]
where $g: \calX \to \bbR^{m}$, the Gauss-Newton/Levenberg-Marquardt methods can be applied, which require evaluation of the vector function $g(x)$ and its Jacobian $J_g(x)$ at each iteration.
In the context of nonlinear least squares problems, the  Levenberg-Marquardt method has shown to be particularly effective, especially when combined with the variable projection
\cite{GolubPereyra73SJoNA-Differentiation,Usevich.Markovsky14A-Optimization}.
\end{itemize}
\end{rem}

The choice of the type of the variable projection depends on a particular problem and its dimensions.
For example, the variable projection with respect to $\widehat{h}$ (Example~\ref{ex:varpro_h}) is reasonable when the degree of the common divisor is small, see  \cite{Karmarkar.Lakshman98JoSC-Approximate,Cheze.etal11TCS-subdivision,Corless.etal95conf-singular,Pan01IaC-Computation,Markovsky.VanHuffel06conf-algorithm,Zhi.etal04JJoIaAM-Hybrid,Li.Zhi13TCS-Computing,Li.etal08TCS-Approximate,Kaltofen.etal08conf-Exact}.
For real polynomials and uniform weights, in  \cite{Markovsky.VanHuffel06conf-algorithm} it was shown that $f_1(\widehat{h})$ and its gradient can be evaluated in $O(dn)$ flops.

Variable projection with respect to $\widehat{\bfg}$ is less common (used in framework of Common Factor Estimation \cite{Stoica.Soderstrom97A-Common,Agrawal.etal04SP-Common}).
But, in fact, the inner minimization problem \eqref{eq:inner_cofe} is well-known in the AGCD literature: it is exactly the so-called least squares division (see, for example, \cite{Zeng11-Numerical}).

\section{Structured low-rank approximation approaches}\label{sec:slra}
In this section, we recall the structured low-rank approximation problem, and review the most popular parameterizations of the problem \eqref{eq:acd2}, adapted to our case.

\subsection{SLRA problem}
The  structured low-rank approximation problem is formulated as follows
\cite{Markovsky.Usevich13SJMAA-Structured}. 
An \textit{affine matrix structure} is an affine map from the \textit{structure
 parameter space} $\bbF^{n_p}$ to the  space of matrices $\bbF^{K \times L}$ (where $\bbF$ is $\bbR$ or $\bbC$), defined by
\begin{equation}
\scrS(p) = S_0 + \sum\limits_{i=1}^{n_p} p_k S_k, \quad S_k \in \bbF^{K\times L}.\end{equation}

\begin{prob}[Structured low-rank approximation]\label{prob:slra1}
Given an affine structure $\scrS$, data vector ${p} \in \bbR^{n_p}$,  and natural number $r < \min(K,L)$ 
\begin{equation}\label{eq:prob_slra}\tag{SLRA}
\minim_{\widehat{p} \in \bbR^{n_p}} \| {p} - \widehat{p} \|_w \;  \sto \;
\rrank \scrS(\widehat{p}) \le r,
\end{equation}
where $\|\cdot\|_w$ is the weighted extended seminorm, see Section~\ref{sec:wnorm}.
\end{prob}

Thus if the are able to represent the set $\pspacegcd{d}$ through the set of low-rank matrices, then Problem~\ref{prob:acd} can be reformulated as Problem~\ref{prob:slra1}.
The classic theorem of Sylvester provides this correspondence for two polynomials.

\begin{theorem}[Sylvester]\label{thm:sylv_2poly}
Two homogeneous polynomials $p \in \pspace{n}$ and $q \in \pspace{m}$ have a non-trivial common divisor if and only if the matrix
\begin{equation}\label{eq:sylv_2poly}
\bmx
\multmat{\pno{p}{2}}{n_1-1} & \multmat{\pno{p}{1}}{n_2-1} 
\emx
\end{equation}
is rank deficient. The rank defect of the matrix is equal to the degree of the GCD.
\end{theorem}

There exists a generalization of Theorem~\ref{thm:sylv_2poly}, for so-called subresultant matrices.

\subsection{Generalized Sylvester subresultant matrix}
For several polynomials (a tuple $\bfp \in \pspace{\bfn}$), the following matrix is typically considered (called generalized Sylvester subresultant matrix) \cite{Rupprecht99JPAA-algorithm,Kaltofen.etal06conf-Approximate}.
\begin{equation}\label{eq:sylv_small}
\sylv{d}(\bfp) \defeq
\bmx
\multmat{\pno{p}{2}}{n_1-d} & -\multmat{\pno{p}{1}}{n_2-d} & 0      & 0 \\
\vdots                      & 0                           & \ddots & 0 \\
\multmat{\pno{p}{N}}{n_1-d} & 0                           & 0      & -\multmat{\pno{p}{1}}{n_N-d}
\emx,
\end{equation}
where $\multmat{p}{k}$ is defined in Section~\ref{sec:pol_mult}. The matrix  $\sylv{d}(\bfp)$  has $K$ rows and $L$ columns:
\begin{equation}\label{eq:sizes_sylv_small}
K = (N-1)(n_1-d+1) + \sum\limits_{k=2}^N n_k, \quad L = \sum\limits_{k=1}^N (n_k-d+1), \quad K \ge L.
\end{equation}
The matrix \eqref{eq:sylv_small} is called the \emph{generalized Sylvester subresultant matrix}. 

\begin{example}\label{ex:sylv_small_3}
For  $\bfp =(\pno{p}{1}, \pno{p}{2}, \pno{p}{3})$, the matrix $\sylv{d}(\bfp)$ has the form
\[
\sylv{d} (\bfp) =
\bmx
\multmat{\pno{p}{2}}{n_1-d}  & -\multmat{\pno{p}{1}}{n_2-d}   &    \bfzero     \\
\multmat{\pno{p}{3}}{n_1-d}  &  \bfzero                   &  -\multmat{\pno{p}{1}}{n_3-d} \\
\emx.
\]

\end{example}

It can be shown that the following lemma holds true. 
\begin{lem}\label{lem:standard_extended_sylvester}
For ${\bfp} = (\pno{p}{1}, \ldots, \pno{p}{N}) \in \pspace{\bfn}$, $\pno{p}{1} \neq \bfzero$, we have that 
\[
{\bfp} \in \pspacegcd{d} \iff\rrank \sylv{d}({\bfp}) \le L - 1,
\]
or, equivalently $\sylv{d}({\bfp})$ has rank deficiency at least $1$.
\end{lem}
\begin{pf}The proof is given in \ref{sec:sylv_full_properties}.\end{pf}

\begin{rem}\label{rem:sylv_subresultant_deficiency}
The following equality of sets holds true:
\[
\{\bfp : \sylv{d}(\bfp)\mbox{ is rank deficient}, \; \pno{p}{1} \neq \bfzero\} = 
\{\bfp \in \pspacegcd{p} : \pno{p}{1} \neq \bfzero\}.
\]
But, the set of $N$-tuples  $\bfp$ with rank-deficient subresultant matrix\begin{equation}
\{ \bfp \in \pspacen \;:\; \sylv{d} (\bfp) \quad \mbox{is rank deficient} \},
\label{eq:ker_repr}
\end{equation}
does not coincide with $\pspacegcd{d}$, in view of our definition of GCD.
Indeed, if $\pno{p}{1} = \bfzero$ then the matrix $\sylv{d}(\bfp)$ becomes automatically rank-deficient.
\end{rem}

\subsection{Full Sylvester subresultant matrix}
By Remark~\ref{rem:sylv_subresultant_deficiency}, if the solution of SLRA f has $\pno{\widehat{p}}{1} = \bfzero$, then it does not give a desired solution to Problem~\ref{pro:acd2}. 
In order to handle properly this non-generic case, we can use an alternative structure, which can be constructed recursively from subresultants.

For an $N$-tuple ${\bfp} = (\pno{p}{1}, \ldots, \pno{p}{N}) \in \pspace{\bfn}$, define
\[
\sylv{d}^{(full)}(\bfp) \defeq
\begin{cases}
\bmx
\sylv{d}(\bfp)\\\hline
\begin{array}{cc}\bfzero & \sylv{d}\big((\pno{p}{2}, \ldots, \pno{p}{N})\big) \end{array} 
\emx, & \mbox{if} \quad N > 2, \\
\sylv{d}(\bfp), & \mbox{if} \quad N = 2,
\end{cases}
\]
where $\sylv{d}(\bfp)$ is defined in \eqref{eq:sylv_small}, and the number of columns of the zero block is $n_1 -d+1$. Thus $\sylv{d}^{(full)}(\bfp)$ has $K^{(full)}$ rows and $L$ columns, where $L$ is as in \eqref{eq:sizes_sylv_small} and
\[
K^{(full)} = \sum\limits_{k=1}^{N-1} \sum\limits_{l=k+1}^N (n_k+n_l-d+1).
\]

\begin{example}[{Example~\ref{ex:sylv_small_3}, continued}]
For  $\bfp =(\pno{p}{1}, \pno{p}{2}, \pno{p}{3})$, we have
\[
\sylv{d}^{(full)} (\bfp) =
\bmx
\multmat{\pno{p}{2}}{n_1-d}  & -\multmat{\pno{p}{1}}{n_2-d}   &    \bfzero     \\
\multmat{\pno{p}{3}}{n_1-d}  &  \bfzero                   &  -\multmat{\pno{p}{1}}{n_3-d} \\
\bfzero  & \multmat{\pno{p}{3}}{n_2-d} & -\multmat{\pno{p}{2}}{n_3-d}        \\
\emx.
\]
\end{example}
\begin{rem}
Compared with the matrix $\sylv{d}(\bfp)$  (which has $N-1$ block rows), the matrix $\sylv{d}^{(full)} (\bfp)$ has $\binom{N}{2}$ block rows, i.e., all possible pairs of the polynomials are present.
The structure of $\sylv{d}^{(full)} (\bfp)$ is similar to the structure of Young flattening of tensors \cite[\S 3.8]{Landsberg12-Tensors}, and probably has a similar algebraic description.
\end{rem}

Next, we show that SLRA of  $\sylv{d}^{(full)}(\bfp)$ is equivalent to Problem~\ref{pro:acd2}.
\begin{prop}\label{prop:sylv_equiv}
For $\bfp = (\pno{p}{1}, \ldots, \pno{p}{N}) \in \pspacen$, we have that
\begin{equation}
 \deg \gcd \bfp \ge d \iff \sylv{d}^{(full)}(\bfp) \mbox{ is rank-deficient}.
 \end{equation}
\end{prop}
\begin{pf}The proof is given in \ref{sec:sylv_full_properties}.\end{pf}

Apart from the precise correspondence between the approximation problems (proved in Proposition~\ref{prop:sylv_equiv}), the structure $\sylv{d}^{(full)}(\bfp)$ can be used to construct initial approximation in the optimization methods (see Section~\ref{sec:iniapprox}).
As shown in the following lemma, the quotient polynomials can be obtained from the kernel of $\sylv{d}^{(full)} (\bfp)$.
\begin{lem}\label{lem:uk_gk}
If $\deg {\rm gcd} ({\bfp}) = d$ and $\sylv{d} ({\bfp}) \bfu = \bfzero$, for $\bfu = (\pno{u}{1},\ldots,\pno{u}{N}) \in \pspace{\bfn-d} \setminus \{\bfzero\}$, then $\pno{p}{k} = \pno{u}{k} h$, where $h \in \gcd({\bfp})$.
\end{lem}
\begin{pf}
The proof is given in \ref{sec:sylv_full_properties}.
\end{pf}


\subsection{Extended Sylvester matrix}
Yet another elegant extension of the Sylvester matrix was proposed in \cite{Karcanias.etal06CMwA-Approximate}. 
The \textit{extended Sylvester matrix}, \cite[(2.2b)]{Karcanias.etal06CMwA-Approximate}) for a parameter $L' \ge \max_k n_k$, defined as
\begin{equation}\label{eq:sylv_alt}
S'_{L'}(\bfp) = \bmx \multmat{\pno{p}{1}}{L'-n_1-1} & \dots  & \multmat{\pno{p}{1}}{L'-n_{N}-1} \emx^{\top}.
\end{equation}
The number of rows $K' = \sum\limits_{k=1}^N (L'-n_k)$ does not exceed the number of columns $L'$.
For the structure \eqref{eq:sylv_alt}, the following theorem holds true.
\begin{theorem}[{\cite[Thm. 1]{Karcanias.etal06CMwA-Approximate}}]\label{thm:sylv_alt}
For  $\bfp \in \pspace{\bfn}$ and $L'\ge  \max_k n_k$,
\begin{equation*}
\rrank S'_{L'}(\bfp) = L' - \deg \gcd \widehat{\bfp}.
\end{equation*}
\end{theorem}
The proof of Theorem~\ref{thm:sylv_alt} (which can be found in \cite[Thm. 1]{Karcanias.etal06CMwA-Approximate}) is based on the following fact on the right kernel of $S'_{L'}(\bfp)$ (an analogue of Lemma~\ref{lem:uk_gk}).

\begin{rem}\label{rem:thm_sylv_alt_pf}
Suppose (for simplicity) that the polynomial $h = \gcd \bfp$ has $d$ simple roots (excluding $\infty$).
From \eqref{eq:sylv_alt},  it follows that for the roots  $\lambda$ of $h$, the vector
\begin{equation}\label{eq:exp_vector}
\bmx 1 & \lambda & \dots & \lambda^{L'-1} \emx^{\top},
\end{equation}
is in the right kernel of $S'_{L'}(\bfp)$. It can be shown that  the vectors in the right kernel are linear combinations of the vectors of the form \eqref{eq:exp_vector},
thus the rank defect of $S'_{L'}(\bfp)$ is $d$.
\end{rem}

Hence, the Problem~\ref{pro:acd2} is equivalent to structured low-rank approximation of the matrix $S'_{L'}(\bfp)$, and can be solved as an SLRA problem.

\begin{rem}
Note that the methods \cite{Karcanias.etal06CMwA-Approximate,Karcanias.etal06IJoC-Matrix} (matrix pencil methodologies) do not solve the structured low-rank approximation problem.
Instead, they are based on unstructured relaxations of the problem (using singular value decomposition). 
\end{rem}

\subsection{Initial approximations for the direct parameterization methods}\label{sec:iniapprox}
The formulation \eqref{eq:prob_slra} also provides a heuristic to obtain an initial guess for $\widehat{h}$ and/or $\widehat{\bfg}$ for the optimization methods in image or kernel representation (from Section~\ref{sec:nls}).
The common heuristic consists in replacing the SLRA problem by unstructured low-rank approximation of a structured matrix $\scrS (\bfp)$.
The unstructured low-rank approximation can be computed, for example, using the SVD. 

The first option is to use  Lemma~\ref{lem:uk_gk}, and use an approximate solution of $\sylv{d} ({\bfp}) \bfu \approx 0$.
Such a solution obtained by unstructured low-rank approximation will be denoted $\bfu_{LRA} \in \pspace{\bfn-d}$.
We may assume that $\bfu_{LRA}$ give an approximation of quotient polynomials (due to the fact that the condition  $\deg {\rm gcd} (\widehat{\bfp}) = d$ represents generic points in $\pspacegcd{n}$).
Then an initial approximation of $\widehat{h}$, can be found by finding the minimizer of \eqref{eq:inner_cofe}.
Let us summarize this option in the following algorithm.
\begin{algorithm}\label{alg:iniapprox_sylv_subr}
Input: $N$-tuple $\bfp \in \pspace{\bfn}$, $d \ge 0$. Output:  initial approximation ${\widehat{h}}_0(z)$.

\begin{itemize}
\item Compute $\bfu_{LRA}$ --- last right singular vector of $\sylv{d}^{(full)}(\bfp)$ (or $\sylv{d}(\bfp)$); 
\item Set ${\widehat{h}}_0 = {\rm arg} \min_{\widehat{h} \in \pspace{d}} \sum\limits_{k=1}^{N} \|\multmat{\phatno{u}{k}}{d} \widehat{h} - \pno{p}{k} \|^2_w$.
\end{itemize}
\end{algorithm}

\begin{rem}
In Algorithm~\ref{alg:iniapprox_sylv_subr}, it may be preferable to use $\sylv{d}^{(full)}(\bfp)$, because it contains  each polynomial $\pno{p}{k}$ the same number of times.
\end{rem}

Another option is to use the approximate kernel of the structure \eqref{eq:sylv_alt}, described in Remark~\ref{rem:thm_sylv_alt_pf},
and compute the initial approximation ${\widehat{h}}_0$ using the matrix pencil approach  (modified matrix pencil method of 
\cite{Karcanias.etal06IJoC-Matrix}). In this case, the algorithm is as follows.

\begin{algorithm}\label{alg:iniapprox_sylv_ext}
Input: $N$-tuple $\bfp \in \pspace{\bfn}$, $d \ge 0$. Output:  initial approximation ${\widehat{h}}_0(z)$.

\begin{itemize}
\item Construct the extended Sylvester matrix $S'_{L'}(\bfp)$ defined in \eqref{eq:sylv_alt}.
\item Compute the SVD $S'_{L'}(\bfp) = U \Sigma V^{*}$, and define by $P \in \bbR^{L' \times d}$ the matrix composed of the  last $d$ columns of $V$.
\item Define $\widehat{Z} = {\rm arg} \min_{Z \in \bbF^{d \times d}} \|P_{1:{L'-1},:} Z - P_{2:L',:}\|_F$.
\item Set  ${\widehat{h}}_0(z) = \det(Z - zI)$ (computed using the eigenvalue decomposition of $Z$).
\end{itemize}
\end{algorithm}

\section{Mosaic Hankel low-rank approximation}\label{sec:mosaic_slra}
In this section, we recall the definition of mosaic-Hankel matrices and results on  mosaic Hankel SLRA \cite{Usevich.Markovsky13JCAM-Variable}. 
Note that compared with \cite{Usevich.Markovsky13JCAM-Variable}, we use transposed matrices.

\subsection{Mosaic Hankel matrices}
Let $\scrH_{k, l}(c) \in \bbF^{k \times l}$ denote a Hankel matrix, generated from  $c \in \bbF^{k+l-1}$, i.e.
\[
\scrH_{k, l}(c) = 
\begin{bmatrix} 
c_1    & c_2     & \cdots    & c_l  \\
c_2    & \iddots & \iddots   &  \vdots \\
\vdots & \iddots & \iddots   & \vdots \\
c_k    & \cdots  & \cdots    & c_{k+l-1}
\end{bmatrix}.
\]
For two vectors $\bfk\in \bbN^{M}$, $\bfl \in \bbN^{T}$, vectors $c^{(i,j)} \in \bbF^{k_i+l_j-1}$, and the combined vector 
\begin{align}
\bfc & = \vcol(c^{(1,1)}, \ldots,c^{(1,T)}, \ldots,c^{(M,1)}, \ldots, c^{(M,T)}) \in \bbF^{n_p}, \nonumber \\
n_p &\defeq \sum\limits_{i,j=1}^{M,T} (k_i+l_j-1),
\label{eq:np_mosaic_hankel}
\end{align}
we define the \textit{mosaic Hankel} matrix \cite{Usevich.Markovsky13JCAM-Variable}:
\begin{equation}\label{eq:mosaic_hankel_structure}
\mathscr{H}_{\bfk, \bfl} (\bfc) \defeq
\bmx
\mathscr{H}_{k_1,l_1} (c^{(1,1)}) & \cdots& \mathscr{H}_{k_1,l_T} (c^{(1,T)})\\
\vdots && \vdots\\
\mathscr{H}_{k_M,l_1} (c^{(M,1)}) & \cdots & \mathscr{H}_{k_M,l_T} (c^{(M,T)})
\emx \in \bbF^{K\times L}.
\end{equation}

\subsection{Structured low-rank approximation and variable projection}
We consider the problem \eqref{eq:prob_slra} for the real-valued structure \eqref{eq:mosaic_hankel_structure} (i.e.,  $\bbF = \bbR$ ). 
We assume that $w \in (0;+\infty]^{n_p}$, $K \ge L$, and denote by $t \defeq L-r$  the rank defect.

Then, following the variable projection approach, as described in \cite{Usevich.Markovsky13JCAM-Variable}, the problem \eqref{eq:prob_slra} can be  rewritten as a bi-level optimization problem:
\begin{align}
&\minim_{P \in \bbR^{L \times t}, \rrank P = t} \optfun{LN}(P),\quad \mbox{where}
\label{eq:slra_kernel}\\
&\optfun{LN}(P) \defeq \left( \min_{\widehat{c} \in \bbR^{n_p}} \|c - \widehat{c}\|^2_{w} \;\sto\;\mathscr{H}_{\bfk, \bfl} (\widehat{c})P = 0 \right).
\label{eq:mult_lln_def_1}
\end{align}
The problem \eqref{eq:mult_lln_def_1} is a linear least norm problem, and has a closed form solution. 
Define $\widehat{c}^{*}(P)$ the optimal solution of \eqref{eq:mult_lln_def_1}, and
\begin{equation}\label{eq:mosaic_slra_lsfunopt}
\lsfunopt{LN}(P) \defeq \diag(\sqrt{w}) (c - \widehat{c}_*(P)),
\end{equation}
such that $\optfun{LN}(P) = \|\lsfunopt{LN}(P)\|^2_2$.
Then the following result holds true.
\begin{theorem}[{\cite[Thm. 1-3]{Usevich.Markovsky13JCAM-Variable}}]\label{thm:complexity}
The complexity (in flops) of the evaluation of $\optfun{LN}(P)$, $\nabla \optfun{LN}(P)$,
$\lsfunopt{LN}(P)$  and the Jacobian of $\lsfunopt{LN}(P)$ with respect to $P$ is 
$O(t^3 L^2 K)$.
\end{theorem}

\begin{rem}
The complexity bound are lower for certain cases (evaluation of $\optfun{LN}(P)$, $\nabla \optfun{LN}(P)$ in the case of uniform weights), but we stick to the bound in Theorem~\ref{thm:complexity}, since it gives complexity for the Gauss-Newton/Levenberg-Marquardt step.
\end{rem}

\begin{rem}
In \eqref{eq:mult_lln_def_1}, the elements of $\lsfunopt{LN}(P)$ corresponding to the infinite weights are equal to $0$, due to equality constraints for elements $\widehat{c}_*$ and $c$ (see also Section~\ref{sec:wnorm}).
\end{rem}

\section{Main results}\label{sec:results}

In this section, we provide the main results of the paper on the connections between the ACD problem and mosaic-Hankel low-rank approximation.

\subsection{Generalized Sylvester LRA as mosaic Hankel LRA}
\label{sec:sylv}
First, we consider the matrix $\sylv{d}(\bfp)$ and show how it can be represented in the form $\mathscr{H}_{\bfk, \bfl} (\bfp) \Phi$, where $\Phi$ is full row rank matrix, and $\mathscr{H}_{\bfk, \bfl} (\bfp)$ is a mosaic Hankel matrix. 
Thus SLRA for this structure can be solved with the methods of \cite{Markovsky.Usevich13JCAM-Software}.

Denote $K$, $L$ as in \eqref{eq:sizes_sylv_small}, $\ell_k \defeq n_k-d$, $\ell_0 \defeq K-n_1-1$, and
\begin{align*}
q^{(1)}  &\defeq \vcol( \bfzero_{\ell_0}, \pno{p}{1}, \bfzero_{\ell_0} ),\quad
q^{(2)}  \defeq \vcol( \bfzero_{\ell_1},\pno{p}{2}, \bfzero_{\ell_1},\pno{p}{3}, \ldots, \bfzero_{\ell_1}, \pno{p}{N}, \bfzero_{\ell_1} ), 
\end{align*} 
where $\bfzero_{n}$ denotes the vector of zeros of length $n$.

\begin{prop}\label{prop:sylv_mosaic}
The  generalized Sylvester subresultant matrix \eqref{eq:sylv_small}  can be represented as the following mosaic Hankel matrix:
\begin{equation}\label{eq:sylv_mosaic}
\sylv{d}(\bfp)= 
\scrH_{K,\bfl} \left(\vcol( q^{(1)} , q^{(2)})\right) \Phi,
\end{equation}
where $\bfl \defeq \bmx \ell_0+1 & \ell_1 + 1\emx^{\top}$ and
\[
\Phi \defeq
\blkdiag \left(\bmx -I_{\ell_N+1} \\ \bfzero_{n_1 \times (\ell_N+1)} \emx, \ldots,
\bmx -I_{\ell_3+1}  \\ \bfzero_{n_1\times (\ell_3+1)} \emx,
-I_{\ell_2+1}, I_{\ell_1+1} \right) J_{L}.
\]
\end{prop}
\begin{pf}
For $p \in \pspace{n}$ and $\ell$ we have
\[
	 \multmat{p}{m}
= \scrH_{\ell+n+1,\ell+1}  \left(\vcol (\bfzero_{\ell},p, \bfzero_{\ell}) \right) J_{\ell+1}.
\]
If we denote $H^{(i,j)} \defeq \scrH_{\ell_i+n_j+1,\ell_i+1}  \left(\vcol (\bfzero_{\ell_i}, \pno{p}{j}, \bfzero_{\ell_i})\right)$,
then we have that
\begin{equation*}
 \sylv{d}(\bfp)=
\bmx
\bfzero & 0 & -H^{(2,1)} & H^{(1,2)}  \\
\bfzero & \iddots & 0  & \vdots\\
-H^{(N,1)} & \bfzero  & 0 & H^{(1,N)}\\
\emx  J= 
\scrH_{K,\bfl} \left(\vcol( q^{(1)} , q^{(2)})\right)  \Phi,
\end{equation*}
which completes the proof.
\end{pf}

\begin{rem}
The  problem \eqref{eq:prob_slra}  for the matrix  \eqref{eq:sylv_small}  can be solved as a weighted mosaic low-rank approximation of the matrix in \eqref{eq:sylv_mosaic}, if we fix the zero elements. This can be accomplished by taking the following weight vector:
\[
w  = 
\vcol(\infty \bfone_{\ell_0}, \pno{w}{1}, \infty \bfone_{\ell_0}, \infty \bfone_{\ell_1}, \pno{w}{2}, \infty \bfone_{\ell_1}, \pno{w}{3}, \ldots, 
\infty \bfone_{\ell_1}, \pno{w}{N}, \infty \bfone_{\ell_1} ),
\]
where $\bfone_{n}$ denotes the vector of ones of length $n$.
We note that the case of infinite weights can be handled by the methods of \cite{Markovsky.Usevich13JCAM-Software}.
\end{rem}

\begin{corollary}\label{cor:sylv_mosaic_complexity}
The complexity of the Gauss-Newton/Levenberg-Marquardt step for the structure \eqref{eq:sylv_mosaic} (using the methods of \cite{Markovsky.Usevich13JCAM-Software}) is
$O(K(K-d+1)^2)$.
\end{corollary}

\begin{rem}\label{rem:varpro_sylv_alt}
As in  Proposition~\ref{prop:sylv_mosaic}, the structure \eqref{eq:sylv_alt} can be represented as a mosaic Hankel structure.
We do not consider this representation here, because SLRA of $S'_{L'}(\bfp)$ presents a difficulty for optimization methods based on kernel representation  of the rank constraint \cite{Markovsky.Usevich13JCAM-Software}, due to the nonlinear structure of the right kernel of a rank deficient $S'_{L'}(\widehat{\bfp})$ (as shown in Remark~\ref{rem:thm_sylv_alt_pf}).
Recently  \cite{Ishteva.etal13-Regularized}, new methods were proposed for structured low-rank approximation of $S'_{L'}(\bfp)$, but the method  \cite{Ishteva.etal13-Regularized}  has cubic computational complexity and is not efficient for large problems. 
\end{rem}

\subsection{Mosaic Hankel matrices and least-squares problems with multiplication matrices}\label{sec:ls_mosaic_relation}
In this section, we establish relations between the variable projection for 
the problem \eqref{eq:agcd_nls} and variable projection for mosaic-Hankel low-rank approximation.

Given  $P \in \bbR^{L \times t}$, an integer vector $\bfk = \rowvec{l_1}{l_T}^{\top} \in \bbN^{T}$ (we also denote $L = \sum\limits_{j=1}^{T} l_j$), and for  $k \ge 1$, we define the matrix-polynomial multiplication matrix:
\begin{equation}
\matmultmat{k}{\bfl} (P) \defeq 
\bmx
\multmat{P^{(1,1)}}{k-1} & \cdots & \multmat{P^{(1,t)}}{k-1} \\
\vdots                   &        & \vdots  \\
\multmat{P^{(T,1)}}{k-1}   & \cdots & \multmat{P^{(T,t)}}{k-1} \\
\emx,
\label{eq:matmultmat}
\end{equation}
where $P^{(i,j)} \in \bbR^{l_i \times 1}$ are the following sub-matrices of the matrix $P$:
\begin{equation}\label{eq:P_partitioning}
P = 
\bmx
P^{(1,1)} & \cdots & P^{(1,t)} \\
\vdots    &        & \vdots    \\
P^{(T,1)} & \cdots & P^{(T,t)}
\emx.
\end{equation}
For two integer vectors $\bfl$ and  $\bfk = \rowvec{k_1}{k_M}^{\top}$, we define 
\[
\matmultmat{\bfk}{\bfl} (P) \defeq \blkdiag \left(\matmultmat{k_1}{\bfl}(P), \ldots, \matmultmat{k_M}{\bfl}(P)\right) \in \bbR^{n_p \times (Kt)}.
\]
where $n_p$ is defined in \eqref{eq:np_mosaic_hankel}, and  $K = \sum\limits_{j=1}^{M} k_j$.
Next,   consider a vector $b \in \bbR^{n_p}$, and a weight vector $v \in [0;+\infty)^{n_p}$. 
Then the following proposition takes place.

\begin{prop}\label{prop:main}
For $w \defeq v^{-1}$ and $c \defeq \diag(v)b$, the solutions of the the minimization problem \eqref{eq:mult_lln_def_1} and the problem 
\begin{equation}
\optfun{LS}(P) \defeq  \min_{x \in \bbR^K} \|\matmultmat{\bfk}{\bfl} (P) x - b\|^2_{v},
\label{eq:mult_lls_def_1}
\end{equation}
are related through the following correspondence:
\begin{equation}\label{eq:optfun_ls_ln}
\optfun{LS}(P) = \|b\|^2_v - \optfun{LN}(P). 
\end{equation}
Moreover, for the  function
\begin{equation}\label{eq:lsfunopt_lls}
\lsfunopt{LS}(P) \defeq \diag{\sqrt{v}} (b - \matmultmat{\bfk}{\bfl} (P) x_*),
\end{equation}
where $x_*$ is the minimizer of \eqref{eq:mult_lls_def_1}, it holds that
\begin{equation}
\lsfunopt{LS}(P) = \diag{\sqrt{v}} b  - \lsfunopt{LN}(P),
\label{eq:lls_lln_sol_2}
\end{equation}
where $\lsfunopt{LN}$ is defined in \eqref{eq:mosaic_slra_lsfunopt}.
Also,  $\optfun{LS}(P) = \|\lsfunopt{LS}(P)\|^2_2$.
\end{prop}
\begin{pf}
From  \cite[Sec.3]{Usevich.Markovsky13JCAM-Variable}), we have that the following equality takes place:
\begin{equation}
\matmultmat{\bfk}{\bfl}^{\top} (P) c = 0 \iff \mathscr{H}_{\bfk, \bfl} (c) P = 0.
\label{eq:mult_mat_low_rank}
\end{equation}
The rest follows from the  correspondence between linear least squares and linear least norm problems presented in \ref{sec:ls_ln}.
\end{pf}

Proposition allows us to apply the results from \cite{Usevich.Markovsky13JCAM-Variable} for the complexity of the Levenberg-Marquardt/Gauss-Newton step in Remark~\ref{rem:algorithms}. (In this case, the function $\lsfunopt{LS}$ from \eqref{eq:lls_sol} can be taken as the function $g$ in Remark~\ref{rem:algorithms}.)

\begin{corollary}
The complexity (in flops) of the evaluation of $\optfun{LS}(P)$, $\nabla \optfun{LS}(P)$,
$\lsfunopt{LS}(P)$  and the Jacobian of $\lsfunopt{LS}(P)$ with respect to $P$ is  $O(t^3 L^2 K)$.
\end{corollary}

\subsection{Variable projection for real-valued polynomials}\label{sec:varpro_real}
In this subsection, we consider the case $\bbF = \bbR$ and the methods presented in Section~\ref{sec:varpro_methods}.
We denote by $n = \sum_k n_k$ the sum of all degrees of the polynomials.

First, we consider variable projection with respect to common divisors (Example~\ref{ex:varpro_h}). In this case, the cost function can be expressed as the cost function \eqref{eq:mult_lls_def_1}
\[
f_1(\widehat{h}) = \optfun{LS} (P),
\]
if we put $P = \widehat{h}$, $x=\bfg$, $\bfl = \bmx d + 1\emx$, $\bfk = \bfn-d$ and $v = \bfw$.
Indeed, in this case,
\[
\matmultmat{\bfk}{\bfl} (P) = \blkdiag(\multmat{\widehat{h}}{n_1 - d}, \ldots,\multmat{\widehat{h}}{n_N - d}).
\]

\begin{corollary}\label{cor:varpro_h_complexity}
The function $f_1$, its gradient, the vector $\lsfunopt{LS}(P)$ from \eqref{eq:lsfunopt_lls}, and its Jacobian can be evaluated in $O(d^2 n)$ flops.
\end{corollary}

Thus the variable projection with respect to $\widehat{h}$ is especially beneficial if $d \ll n_k$.

\begin{rem}
Note that for uniform weights (when $\pno{w}{k} = const$ in \eqref{eq:wdist}), the complexity is $O(dn)$ (which corresponds to the results of \cite{Corless.etal95conf-singular,Markovsky.VanHuffel06conf-algorithm}). 
\end{rem}

\begin{rem}\label{rem:varpro_h_invariance}
The cost function $f_1$ is invariant to multiplication by a scalar, i.e. $f_1(\alpha \widehat{h}) = f_1(\widehat{h})$  for any $\alpha \neq 0$, since the columns of $\matmultmat{\bfk}{\bfl} (\widehat{h})$ and $\matmultmat{\bfk}{\bfl} (\alpha \widehat{h})$ span the same subspace. Therefore, \eqref{eq:outer} is a minimization on the projective space (which is a special case of minimization  the Grassmann manifold \cite{Usevich.Markovsky14A-Optimization}). 

Moreover, due to correspondence \eqref{eq:optfun_ls_ln}, the minimization problem \eqref{eq:outer} has a  similar structure to \eqref{eq:slra_kernel}.
 This fact is exploited in the software package \cite{Markovsky.Usevich13JCAM-Software}. 
\end{rem}

Now let us consider variable projection with respect to $\bfg$  (Example~\ref{ex:varpro_bfg}). In this case, the cost function can be expressed as 
\[
f_2(\widehat{\bfg}) = \optfun{LS} (P),
\]
where $P = \widehat{\bfg}$,  $\bfk = \bmx d + 1\emx$ and $\bfl = \bfn-d$ and $v = \bfw$.

\begin{corollary}\label{cor:cofe_complexity}
The function $f_2$, its gradient, the vector $\lsfunopt{LS}(P)$ from \eqref{eq:lsfunopt_lls}, and its Jacobian can be evaluated in $O((n-Nd)^2 n)$. 
\end{corollary}
The variable projection with respect to $\bfg$ is beneficial if $(n_k - d) \ll n$ (i.e., where the degrees of the quotients are fixed and small). 
In this case, the complexity is linear in the degrees of the polynomials.
\begin{rem}
For $N > 2$, the complexity in Corollary~\ref{cor:cofe_complexity} is smaller than in the case of the kernel representation approach (using generalized Sylvester subresultant matrices), where the complexity never approaches the linear rate, see the Corollary~\ref{cor:sylv_mosaic_complexity}.
\end{rem}

\begin{rem}
The cost function, as in Remark~\ref{rem:varpro_h_invariance} is invariant with respect to scaling, i.e., $f_2(\alpha \widehat{\bfg}) = f_2(\widehat{\bfg})$ for $\alpha \neq 0$, since the columns of $\matmultmat{\bfk}{\bfl} (\widehat{\bfg})$ and $\matmultmat{\bfk}{\bfl} (\alpha \widehat{\bfg})$ span the same subspace. 
Therefore, problem \eqref{eq:outer_cofe} has similar structure to problem~\eqref{eq:slra_kernel}, which can be solved by the software package \cite{Markovsky.Usevich13JCAM-Software}. 
\end{rem}

\subsection{Variable projection for complex-valued polynomials}\label{sec:varpro_complex}
Now  assume that the polynomials in \eqref{eq:agcd_nls} are complex. (In this section, we only consider variable projection w.r.t. $\widehat{h}$.) Then the polynomials can be represented as 
\[
\pno{p}{k} = \pno{p}{\mathscr{R},k} + i \cdot \pno{p}{\mathscr{I},k}, \quad
\widehat{h} = \widehat{h}^{\mathscr{R}} + i \cdot \widehat{h}^{\mathscr{I}},\quad
\pno{g}{k} = \pno{g}{\mathscr{R},k} + i \cdot \pno{g}{\mathscr{I},k},
\]
where $i = \sqrt{-1}$. If we set 
\[
\begin{split}
&
\bfp^{(\bbR)} = \vcol \left(\pno{p}{\mathscr{R},1}, \pno{p}{\mathscr{I},1}, \ldots, \pno{p}{\mathscr{R},N}, \pno{p}{\mathscr{I},N} \right), \\
& \bfg^{(\bbR)} = \vcol \left(\pno{g}{\mathscr{R},1}, \pno{g}{\mathscr{I},1}, \ldots, \pno{g}{\mathscr{R},N}, \pno{g}{\mathscr{I},N} \right), \\
& \bfw^{(\bbR)} = \vcol \left(\pno{w}{1}, \pno{w}{1}, \ldots, \pno{w}{N}, \pno{w}{N} \right), 
\end{split}
\]
\[
\begin{split}
& P^{(\widehat{h})} \defeq 
\bmx
\widehat{h}^{\mathscr{R}} & -\widehat{h}^{\mathscr{I}} \\
\widehat{h}^{\mathscr{I}} & \widehat{h}^{\mathscr{R}} \\
\emx \in \bbR^{(d+1)\times 2}.
\end{split}
\]
Then we have that the problem \eqref{eq:agcd_nls} is  equivalent to
\begin{equation}
\min_{\widehat{h},\widehat{\bfg}} \|\matmultmat{\bfk}{\bfl} (P^{(\widehat{h})}) \bfp^{(\bbR)} - \bfg^{(\bbR)}\|^2_{\bfw^{(\bbR)}},
\label{eq:nls_complex}
\end{equation}
for $\bfl = \bmx d + 1 & d+1\emx^{\top}, \bfk = \bfn-d, v = \bfw$.
Then \eqref{eq:nls_complex} can be  rewritten as
\begin{align}
&\minim_{\widehat{h} \in \pspace{d} \setminus \{0\}} f_3(P^{(\widehat{h})}), \quad \mbox{where}  \label{eq:outer_c} \\
&f_3(P^{(\widehat{h})}) \defeq \mbox{minimum of \eqref{eq:nls_complex} for fixed}\;P^{(\widehat{h})} = \optfun{LS} (P^{(\widehat{h})}).
\end{align}
It can be seen that $f_3(P^{(\widehat{h})}) = f_3(P^{(\alpha \widehat{h})})$ for any $\alpha \in \bbC \setminus \{0\}$. Therefore, $f_3$ can be minimized on the complex projective plane.
Another option is to use parameterization
\[
\mvec\left((P^{(\widehat{h})})^{\top}\right) = 
\bmx
A & -B \\ B & A
\emx \bmx
\widehat{h}^{\mathscr{R}} \\
\widehat{h}^{\mathscr{I}}  \\
\emx,
\quad A \defeq I_{d+1} \otimes \bmx 1\\0\emx, \quad B \defeq I_{d+1} \otimes \bmx 0\\1\emx,
\]
which is supported in the  package \cite{Markovsky.Usevich13JCAM-Software}, and we  use it in numerical experiments.

\subsection{Accuracy of the computations}\label{sec:cond}
As shown in \cite{Usevich.Markovsky13JCAM-Variable}, the key step in evaluation of \eqref{eq:optfun_ls_ln}, \eqref{eq:lls_lln_sol_2} and their derivatives is the solution of a system of equations  $\Gamma u = v$, where  $\Gamma$ is defined in \cite{Usevich.Markovsky13JCAM-Variable}. From  Proposition~\ref{prop:main}, we have that 
\[
\Gamma (P) = \blkdiag\left(\Gamma^{(k_1)} (P), \ldots, \Gamma^{(k_M)} (P)\right),
\]
where $\Gamma^{(k)} (P) \defeq \matmultmat{k}{\bfl}^{\top} (P)\diag{v}  \matmultmat{k}{\bfl} (P)$. In the software \cite{SLRA},  the system $\Gamma u = v$ is solved block by block, using Cholesky factorization. 
The accuracy of solving the subsystems $\Gamma^{(k)} (P) u = v$  mainly depends on the condition number of $\Gamma^{(k)}$ \cite{Golub.VanLoan96-Matrix}. 
\begin{rem}
Although in this paper we use Cholesky factorisation, the QR factorisation (for example, using the updating strategy of \cite{Zeng11-Numerical}) may be used to avoid squaring the condition number.
\end{rem}

In what follows, for simplicity, we consider the case of $2$-norm ($v \equiv 1$). (In fact, the case of blockwise weighted $2$-norm is similar \cite{Usevich.Markovsky13JCAM-Variable}.) In this case, the matrix  $\Gamma^{(k)} (P) $ is block-Toeplitz, and  behaviour of its eigenvalues depends on its \textit{symbol} \cite{Miranda.Tilli00SJMAA-Asymptotic}
\begin{equation}
\rmF(z) = P^{*} (z)  P (z),
\label{eq:f_cond}
\end{equation}
where $P(z)$ is  the matrix polynomial 
\[
P(z) \defeq 
\bmx
P^{(1,1)} (z) & \cdots & P^{(1,t)} (z) \\
\vdots                 &        & \vdots  \\
P^{(K,1)} (z) & \cdots & P^{(K,t)} (z) \\
\emx,
\]
where $P^{(i,j)} (z)$ are the generating functions of the vectors $P^{(i,j)}$ in \eqref{eq:P_partitioning}. 
Since $\rmF(z)$ is Hermitian for all $z$, and is continuous on the unit circle $\bbT$, we can define
\[
a_\rmF \defeq \min_{\bbT} \lambda_{min} \left(\rmF(z)\right) \ge 0, \quad
b_\rmF \defeq \max_{\bbT} \lambda_{max} \left(\rmF(z)\right) < \infty.
\]
Where $\lambda_{min}(B)$ and $\lambda_{max}(B)$ are the minimal and maximal eigenvalues of a matrix $B$. The results of \cite{Miranda.Tilli00SJMAA-Asymptotic} imply that 
$\lambda_{min} (\Gamma^{(k)}) \searrow a_{\rmF} \quad \mbox{and} \quad
\lambda_{max} (\Gamma^{(k)}) \nearrow b_{\rmF}$,
\ie~the eigenvalues are in the interval $[a_\rmF;b_\rmF]$ and converge to the endpoints as $k \rightarrow \infty$. Therefore, the condition number $\kappa_2(\Gamma^{(l)}) \defeq \lambda_{max} (\Gamma^{(k)}) / \lambda_{min} (\Gamma^{(k)})$ behaves as
\begin{equation}\label{eq:condition_est}
\kappa_2(\Gamma^{(k)}) \nearrow \frac{b_{\rmF}}{a_\rmF}.
\end{equation}
If $\rmF(z)$ is positive definite on $\bbT$, then $\kappa_{2}(\Gamma^{(k)}) \le {b_{\rmF}}/{a_\rmF} < \infty$. Otherwise, $\kappa_{2}(\Gamma^{(k)}) \rightarrow \infty$ (results on the rate of convergence are known).
See \cite[\S 6.2]{Usevich.Markovsky13JCAM-Variable} for more details.

\begin{example}\label{ex:cond1}
In the real-valued  case (see Section~\ref{sec:varpro_real}), for variable projection with respect to $\widehat{h}$ (see Example~\ref{ex:varpro_h}), easy calculations show that \eqref{eq:f_cond} becomes
\[
\rmF(z) = |\widehat{h}(z)|^2, \quad \mbox{for } z \in \bbT,
\]
and 
\begin{equation}\label{eq:condition_est_real}
a_{\rmF} = \min_{z \in \bbT} |\widehat{h}(z)|^2, \quad b_{\rmF} = \max_{z \in \bbT} |\widehat{h}(z)|^2.
\end{equation}
From \eqref{eq:condition_est}, we conclude that the computations are well-conditioned if the tentative common divisor $\widehat{h}$ (during optimisation) does not have roots on the unit circle. If it has roots on the unit circle, the computations may become ill-conditioned.
\end{example}
The results in Example~\ref{ex:cond1} are in agreement with the similar analysis of conditioning of $\Gamma^{(k)}$, performed in \cite{Corless.etal95conf-singular}.
Easy calculations show that   \eqref{eq:condition_est}  with \eqref{eq:condition_est_real} are  valid also for complex polynomials (see Section~\ref{sec:varpro_complex})).

\begin{example}\label{ex:cond2}
In the real-valued  case (see Section~\ref{sec:varpro_real}), for variable projection with respect to $\widehat{\bfg}$ (see Example~\ref{ex:varpro_bfg}), \eqref{eq:f_cond} becomes
\[
F(z) = \sum\limits_{j=1}^{N} |\widehat{g}^{(j)}(z)|^2, \quad \mbox{for } z \in \bbT,
\]
and 
\[
a_{\rmF} = \min_{z \in \bbT} \rmF(z), \quad b_{\rmF} = \max_{z \in \bbT} \rmF(z).
\]
Thus the computations in this case may become ill conditioned if all the tentative quotient polynomials $\widehat{g}^{(j)}(z)$ have common roots on the unit circle, which is less likely.
\end{example}

\subsection{Angles between polynomials as an approximation criterion}\label{sec:angles}
Finally, the variable projection principle helps us to understand the importance of the relative distance, which was used by many authors, see for example the discussion in  \cite[\S 4]{Bini.Boito10-fast}.
Define two distances between tuples of polynomials
\begin{equation}
\displaystyle \mathrm{dist}_{sin} (\bfp,\bfq) \defeq \sum\limits_{k=1}^N \sin^2 (\angle(\pno{p}{k}, \pno{q}{k})).
\label{eq:dist_sin}
\end{equation}
and
\begin{equation}
\displaystyle \mathrm{dist}_{nrm} (\bfp,\bfq) \defeq \frac{\|\pno{p}{k}- \pno{{q}}{k}\|_2^2}{\|\pno{p}{k}\|_2^2}.\label{eq:dist_nrm}
\end{equation}
The distance \eqref{eq:dist_sin} does not depend on the scaling of coefficients of $\pno{p}{k}$ and $\pno{q}{k}$ and depends only on the roots of the polynomials. But,  the distance  $\mathrm{dist}_{sin}$  may be difficult to minimize as is.
However, the normalized distance \eqref{eq:dist_nrm} is just a special case of the weighted $2$-norm, but depends on scaling of the polynomials.

In what follows, we prove that the solutions of the problem \eqref{eq:agcd_nls} coincide  for the two  distances. 
The proof is quite simple, but we could not find it in the AGCD literature.

\begin{prop}
For any tuple of polynomials $\bfp$ and for any $d$ we have that
\[
\min_{\footnotesize\begin{array}{c}\widehat{\bfg} \in \pspace{\bfn-d}, \\
\widehat{h} \in \pspace{d}\end{array}}  \mathrm{dist}_{nrm} (\bfp,\widehat{\bfg}\cdot\widehat{h}) = \min_{\footnotesize\begin{array}{c}\widehat{\bfg} \in \pspace{\bfn-d}, \\
\widehat{h} \in \pspace{d}\end{array}}  \mathrm{dist}_{sin} (\bfp,\widehat{\bfg}\cdot\widehat{h}).
\]
\end{prop}
\begin{pf}
Consider minimization of \eqref{eq:dist_sin}. 
By applying the variable projection principle, the problem becomes
\[
f_4(\widehat{h}) := \min_{\widehat{\bfg}} \sum\limits_{k=1}^N \sin^2 (\angle(\pno{p}{k}, \widehat{h}\cdot \pno{\widehat{g}}{k} )) = 
\sum\limits_{k=1}^N \min_{\pno{\widehat{g}}{k}} \sin^2 (\angle(\pno{p}{k}, \multmat{\widehat{h}}{n_k-d} \pno{\widehat{g}}{k})).
\]
It is well known \cite[\S 17.26]{Kendall.Stuart77v2-advanced}, that the least-squares solution minimizes the angle between the approximating vector and the given vector. Let $\pno{\widehat{g}}{k}_{*}$ be the solution of a least squares problem with matrix $A_k = \multmat{\widehat{h}}{n_k-d}$ and right-hand side $\pno{p}{k}$. Then 
$(\pno{p}{k}- A_k \pno{\widehat{g}}{k}_{*}) \bot A_k \pno{\widehat{g}}{k}_{*}$ and we have that
\[
f_4(\widehat{h}) := 
\sum\limits_{k=1}^N
\frac{\| A_k \pno{\widehat{g}}{k}_{*} -  \pno{p}{k} \|_2^2}{\|\pno{p}{k}\|_2^2} = 
\min_{\widehat{\bfg}}  \mathrm{dist}_{nrm} (\bfp,\widehat{\bfg}\cdot\widehat{h}),
\]
which completes the proof.
\end{pf}

\begin{rem}
Minimizing the relative distance \eqref{eq:dist_nrm}  is equivalent to minimizing the $\ell_2$-norm after a preliminary scaling of the input polynomials.
\end{rem}

\section{Numerical examples}\label{sec:exper}

In this section, we provide numerical experiments that include comparison with the state-of-the-art methods.
The methods developed in this paper are implemented in MATLAB and  are based on the SLRA package \cite{SLRA} described in \cite{Markovsky.Usevich13JCAM-Software}. The source code of the methods and experiments is publicly available at \url{http://github.com/slra/slra}. 

In the experiments, the method used for minimization of  $\optfun{LN}(P) = \|\lsfunopt{LN}(P)\|^2_2$ or  $\optfun{LS}(P) = \|\lsfunopt{LS}(P)\|^2_2$ is the Levenberg-Marquardt method. 
In the SLRA package \cite{SLRA}, two implementations of the Levenberg-Marquardt method are currently used: a standard implementation in GNU Scientific Library \cite{GSL} and own implementation that uses data-driven local coordinates approach \cite{Usevich.Markovsky14A-Optimization}, based on the variant described in \cite[p.366]{Pintelon.Schoukens12-System}. In this paper, the former variant  is mainly used in real-valued case, and the latter in complex-valued case.

\subsection{Example of ill-conditioned polynomials}
First, we consider a classic example from \cite[Test 2]{Zeng11-Numerical}  (which can be also found  in \cite[Example 4.2]{Bini.Boito10-fast} and   \cite[Test~5]{Terui13TCS-GPGCD}). The following two polynomials are considered:
\[
u (x) = \prod\limits_{j=1}^{10}(x-x_j), \quad  v (x) = \prod\limits_{j=1}^{10}(x-x_j+10^{-j}), \quad x_j = (-1)^{j} (j/2).
\]
Compared with the mentioned reference, we normalize the polynomials as
\[
p^{(1)} = u(x) / \|u\|_2, \quad p^{(2)} = v(x) / \|u\|_2,
\]
and compare methods according to Euclidean distance $\|(\pno{p}{1}, \pno{p}{2}) - (\pno{\widehat{p}}{1}_*,\pno{\widehat{p}}{2}_*)\|_2$.

All the methods are started from the same initial approximation computed in Algorithm~\ref{alg:iniapprox_sylv_subr}.  ``LRA'' stands for using initial approximation in Section~\ref{sec:iniapprox} without optimization (refinement).  ``VP$_{{h}}$'' denotes the variable projection method w.r.t. $\widehat{h}$ in the image representation (Example~\ref{ex:varpro_h}). ``VP$_{g}$'' denotes the variable projection method w.r.t. $\widehat{\bfg}$  (Example~\ref{ex:varpro_bfg}). ``VP$_{\calS}$''' stands for the variable projection method in the kernel representation (Section~\ref{sec:sylv}).  ``FASTGCD'' denotes the combination of the Gauss-Newton method and line search, used in \cite{Bini.Boito10-fast} (function {\tt c\_f\_newton\_iter}).

\begin{table}
\caption{Optimal approximations of the methods}
\label{tab:ex11}
\centering
{\small
\begin{filecontents*}{res_ex1_matlab.txt}
		d	LRA	VP$_{h}$	VP$_{g}$	VP$_{S}$	UVGCD	FASTGCD
   1.0000000e+00   2.3907363e-03   3.8329847e-16   1.6577868e-15   1.7395454e-16   3.4523550e-16   9.0602906e-16
   2.0000000e+00   1.1003688e-04   2.2485205e-14   7.8135464e-14   8.3530996e-17   2.2483688e-14   8.2953021e-14
   3.0000000e+00   1.0712296e-04   1.5252361e-12   1.5252362e-12   1.8876537e-16   1.5252411e-12   3.0778226e-12
   4.0000000e+00   1.7588879e-05   8.3986866e-11   8.3986867e-11   1.5523331e-14   8.3986869e-11   1.7899372e-10
   5.0000000e+00   4.5847953e-05   4.4865473e-09   4.4865473e-09   1.7209732e-05   4.4865472e-09   7.6870364e-09
   6.0000000e+00   3.2500008e-05   1.8292962e-07   1.8292962e-07   2.7042724e-05   1.8292962e-07   3.3897536e-07
   7.0000000e+00   4.8913826e-04   7.0890249e-06   7.0890249e-06   3.9280758e-04   7.0890249e-06   1.1579774e-05
   8.0000000e+00   2.8472625e-03   2.1037972e-04   1.7288122e-04   1.7288122e-04   1.7288122e-04   2.5507774e-04
   9.0000000e+00   2.1784493e-02   3.9963886e-03   3.9963886e-03   3.9963886e-03   3.9963888e-03   4.4374487e-03
   1.0000000e+01   6.5734913e-02   6.5734913e-02   6.5734913e-02   6.5734913e-02   6.5734913e-02   1.4429767e-01
\end{filecontents*}
\pgfplotstabletypeset[
every head row/.style={after row=\hline},
every first column/.style={column type/.add={}{|}}]{res_ex1_matlab.txt}
}
\end{table}

\begin{table}
\centering
\parbox{.45\linewidth}{
\caption{Number of iterations of the methods}
\label{tab:ex1_iter}
\centering
{\small
\begin{filecontents*}{iters_ex1_matlab.txt}
		d	LRA	VP$_{h}$	VP$_{g}$	VP$_{S}$	UVGCD	FASTGCD
   1.0000000e+00  -1.0000000e+00   4.0000000e+00   5.0000000e+00   3.0000000e+00  -1.0000000e+00   4.0000000e+00
   2.0000000e+00  -1.0000000e+00   4.0000000e+00   4.0000000e+00   2.0000000e+00  -1.0000000e+00   7.0000000e+00
   3.0000000e+00  -1.0000000e+00   4.0000000e+00   4.0000000e+00   2.0000000e+00  -1.0000000e+00   1.3000000e+01
   4.0000000e+00  -1.0000000e+00   4.0000000e+00   4.0000000e+00   9.0000000e+00  -1.0000000e+00   5.0000000e+00
   5.0000000e+00  -1.0000000e+00   5.0000000e+00   4.0000000e+00   6.0000000e+00  -1.0000000e+00   1.2000000e+01
   6.0000000e+00  -1.0000000e+00   6.0000000e+00   4.0000000e+00   1.6000000e+01  -1.0000000e+00   6.0000000e+00
   7.0000000e+00  -1.0000000e+00   6.0000000e+00   4.0000000e+00   1.0000000e+02  -1.0000000e+00   4.0000000e+00
   8.0000000e+00  -1.0000000e+00   3.1000000e+01   3.0000000e+00   5.0000000e+00  -1.0000000e+00   7.0000000e+00
   9.0000000e+00  -1.0000000e+00   4.0000000e+00   3.0000000e+00   4.0000000e+00  -1.0000000e+00   1.9000000e+01
   1.0000000e+01  -1.0000000e+00   1.0000000e+00   1.0000000e+00   1.0000000e+00  -1.0000000e+00   2.0000000e+00
\end{filecontents*}
\pgfplotstabletypeset[
every head row/.style={after row=\hline},
columns={{d},{VP$_{h}$},{VP$_{g}$},{VP$_{S}$},FASTGCD},
every first column/.style={column type/.add={}{|}}]{iters_ex1_matlab.txt}}
}
\quad
\parbox{.45\linewidth}{
\caption{Condition number of the $\Gamma$ matrices}
\label{tab:ex1_cond}
\centering
{\small
\begin{filecontents*}{conds_ex1_matlab.txt}
	d	VP$_{h}$	VP$_{g}$	VP$_{S}$
   1.0000000e+00   2.1697517e+00   1.1896762e+00   7.8855593e+33
   2.0000000e+00   1.2169992e+00   4.6229224e+00   7.4435985e+28
   3.0000000e+00   2.3727060e+00   6.0957080e+00   1.4324572e+24
   4.0000000e+00   1.6080280e+00   1.0126917e+01   2.1470127e+19
   5.0000000e+00   2.9745191e+00   1.3972194e+01   1.6724308e+15
   6.0000000e+00   2.3557235e+00   1.4020119e+01   8.3411812e+10
   7.0000000e+00   3.8195842e+00   1.7061415e+01   7.0685667e+07
   8.0000000e+00   2.9918360e+00   1.2329925e+01   4.9153705e+04
   9.0000000e+00   1.4086507e+00   9.3731983e+00   8.5769756e+02
   1.0000000e+01   1.0000000e+00   1.0000000e+00   1.0000000e+00
\end{filecontents*}
\pgfplotstabletypeset[
every head row/.style={after row=\hline},
every first column/.style={column type/.add={}{|}}]{conds_ex1_matlab.txt}}
}
\end{table}

The results in Table~\ref{tab:ex11} show the Euclidean distances, and in Table~\ref{tab:ex1_iter}, we show number iterations of the methods (unfortunately, the number of iterations for ``UVGCD'' is not available).
In Table~\ref{tab:ex1_cond}, we present the condition numbers for the $\Gamma$ matrix  in the variable projection methods.

The results of the experiments show ``UVGCD'' gives the overall best approximation, and we use its results as a reference.
We see that the methods ``VP$_{{h}}$'' and ``VP$_{{g}}$'' match the results of ``UVGCD'' in few iterations, except the cases $d = 8$ for ``VP$_{{h}}$'' and $d=1,2$ for ``VP$_{{g}}$''.
However, each of these  ``bad'' cases corresponds to the cases where the methods  large search space and high computational complexity, and should not normally be used. 
It is natural to use ``VP$_{h}$'' for $d < \frac{n_k}{2}$ and ``VP$_{h}$'' for $d \ge \frac{n_k}{2}$, (for example, ``VP$_{{g}}$'' in the case $d=8$). 
The $\Gamma$ matrices (see Section~\ref{sec:cond}) are well-conditioned in all cases for these methods.

The method ``FASTGCD'' does not match the results of ``UVGCD'', probably due to the settings of the stopping criteria. Also, the method produces polynomials with complex coefficients as a result.

The method based on the kernel representation (Section~\ref{sec:sylv}) fails to produce good results for $d < 8$. For $d=1,\ldots,4$ it seems to give a good result, but the $\Gamma$ matrix used in computations is ill-conditioned, and the regularization of $\Gamma$ in the package \cite{SLRA} is automatically applied in this example. (This means that the computed approximating polynomials are not guaranteed to have a common divisor.) Therefore, the variable projection in kernel representation should be used only for small $d$.


\subsection{Complex polynomials and speed of the computations}
In this section, we compare speed  of methods ``VP$_{ h}$'' , ``VP$_{g}$'', and  ``FASTGCD''. 

\subsubsection{Small GCD degree scenario}
We consider the example of complex polynomials  from \cite[\S~4.6]{Bini.Boito10-fast}:
\[
\begin{split}
& \pno{q}{1,k}(z) := h(z) \pno{g}{1,k}(z),\quad \pno{q}{2,k} := h(z) \pno{g}{2,k} (z), \\
& h := z^4 + 10z^2 + z - 1, \\
& \pno{g}{1,k}(z) := (z^{k\ell_1} -1) (z^{k\ell_2} -2) (z^{k\ell_3} -3), \\
& \pno{g}{2,k}(z) := (z^{k\ell_1} +i) (z^{k\ell_2} +5) (z^{k\ell_3} +2 ), \\
\end{split}
\]
where $\ell_1 = 25$, $\ell_2 = 15$, $\ell_3 = 10$. We compute normalized polynomials
($ \pno{\widetilde{q}}{1,k}$ and $ \pno{\widetilde{q}}{2,k}$), and also add a small noise to the polynomials:
\[
\pno{p}{1,k} := \pno{q}{1,k} + \varepsilon_{1,k}, \quad \pno{p}{2,k} := \pno{q}{2,k} + \varepsilon_{2,k},
\]
where each $\varepsilon_{j,k}$ is a  realization of the Gaussian zero-mean i.i.d. random vector with standard deviation $10^{-4}$.
We average the results over $20$ realizations of noise.

We consider the test polynomials for $k = \{1,\ldots,8\}$, thus the degrees of the polynomials range between $50$ and $400$. We compare two methods: ``VP$_{ h}$''  and ``FASTGCD''. All the methods are started from the same initial approximation.

\begin{table}
\parbox{.60\linewidth}{
\caption{Optimal approximations of the methods}
\label{tab:res_ex2}
\centering
{\small
\begin{filecontents*}{res_ex2_matlab.txt}
		k	n	LRA	VP$_{h}$(2)	VP$_{h}$	FASTGCD	FASTGCD(1)
   1.0000000e+00   5.5000000e+01   8.5893202e-03   1.6414999e-04   1.6414999e-04   1.5219957e-03   1.6148284e-03
   2.0000000e+00   1.0500000e+02   9.8222717e-03   1.8041256e-04   1.8041256e-04   1.5124317e-03   1.5961913e-03
   3.0000000e+00   1.5500000e+02   7.7303064e-03   1.5874043e-04   1.5874043e-04   3.1522411e-03   4.1694736e-03
   4.0000000e+00   2.0500000e+02   8.5125103e-03   1.5426189e-04   1.5426189e-04   1.9131652e-03   1.9462132e-03
   5.0000000e+00   2.5500000e+02   8.6711632e-03   1.6052464e-04   1.6052464e-04   2.0995456e-03   2.4047724e-03
   6.0000000e+00   3.0500000e+02   7.5545170e-03   1.5694492e-04   1.5694492e-04   2.2099639e-03   2.3186850e-03
   7.0000000e+00   3.5500000e+02   1.1653711e-02   1.6828563e-04   1.6828563e-04   3.4370680e-03   3.7441911e-03
   8.0000000e+00   4.0500000e+02   1.0086329e-02   1.5844044e-04   1.5844044e-04   3.6716157e-03   4.1872957e-03
\end{filecontents*}
\pgfplotstabletypeset[
every head row/.style={after row=\hline},
columns={{k},{n},{LRA},{VP$_{h}$},FASTGCD},
every first column/.style={column type/.add={}{|}}]{res_ex2_matlab.txt}
}}
\quad
\parbox{.35\linewidth}{
\caption{Average number of iterations}
\label{tab:iters_ex2}
\centering
{\small
\begin{filecontents*}{iters_ex2_matlab.txt}
		k	n	LRA	VP$_{h}$(2)	VP$_{h}$	FASTGCD	FASTGCD(1)
   1.0000000e+00   5.5000000e+01  -1.0000000e+00   2.0000000e+00   2.1000000e+00   3.3000000e+00   1.0000000e+00
   2.0000000e+00   1.0500000e+02  -1.0000000e+00   2.0000000e+00   2.1000000e+00   3.6000000e+00   1.0000000e+00
   3.0000000e+00   1.5500000e+02  -1.0000000e+00   2.0000000e+00   2.1000000e+00   3.3500000e+00   1.0000000e+00
   4.0000000e+00   2.0500000e+02  -1.0000000e+00   2.0000000e+00   2.0500000e+00   3.3500000e+00   1.0000000e+00
   5.0000000e+00   2.5500000e+02  -1.0000000e+00   2.0000000e+00   2.0500000e+00   3.3000000e+00   1.0000000e+00
   6.0000000e+00   3.0500000e+02  -1.0000000e+00   2.0000000e+00   2.0000000e+00   2.9500000e+00   1.0000000e+00
   7.0000000e+00   3.5500000e+02  -1.0000000e+00   2.0000000e+00   2.2500000e+00   2.9500000e+00   1.0000000e+00
   8.0000000e+00   4.0500000e+02  -1.0000000e+00   2.0000000e+00   2.0500000e+00   3.6000000e+00   1.0000000e+00
\end{filecontents*}
\pgfplotstabletypeset[
every head row/.style={after row=\hline},
columns={{k},{n},{VP$_{h}$},FASTGCD},
every first column/.style={column type/.add={}{|}}]{iters_ex2_matlab.txt}
}}
\end{table}

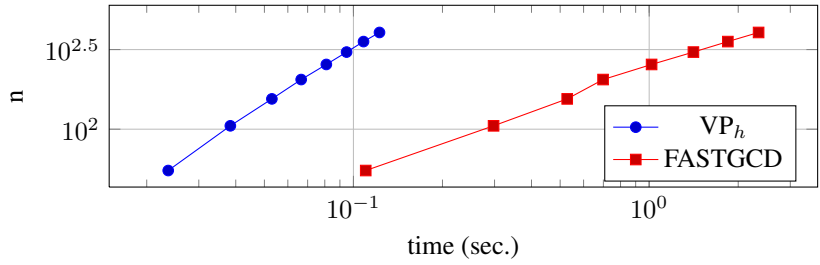
\begin{figure}[ht!]
\centering
\begin{tikzpicture}
\begin{loglogaxis}[ylabel=n,xlabel=time (sec.),ymax=600,legend pos = south east, height = 4cm, width = 11cm, grid=major]
\pgfplotstableread{times_ex2_matlab.txt}\loadedtable %
\addplot table[y index=1,x index=2] {\loadedtable};
\addplot table[y index=1,x index=3] {\loadedtable};
\legend{VP$_{{h}}$,FASTGCD}
\end{loglogaxis}
\end{tikzpicture}
\caption{Comparison of the time per iteration, small $d$} 
\label{fig:speed2}                                 
\end{figure}

As shown in Tables~\ref{tab:res_ex2} and~\ref{tab:iters_ex2}, the method ``VP$_{ h}$'' achieves better  approximation error with similar number of iterations.  
For measuring speed, we limit the number of iterations to $1$ in ``VP$_{ h}$''   and call directly function {\tt c\_iterfast} (one iteration of ``FASTGCD'').   In Fig.~\ref{fig:speed2}, the time  is plotted versus~$n$.

\subsubsection{Large GCD degree scenario}
We repeat the same experiments, but for the case of growing GCD degree. We consider polynomials from  \cite[Ex.~4.3]{Bini.Boito10-fast}:
\[
\begin{split}
& \pno{q}{1}(z) := (\sum\limits_{j=0}^3 z^j)u_d(z),\quad \pno{q}{2}(z) := (\sum\limits_{j=0}^3 (-z)^j)u_d(z),
\end{split}
\]
and $u_d$ is a polynomial of degree $d$, whose coefficients are random integers in $[-5;5]$. We compute normalized polynomials
($ \pno{\widetilde{q}}{1}$ and $ \pno{\widetilde{q}}{2}$), and add a small noise:
\[
\pno{p}{1} := \pno{q}{1} + \varepsilon_{1}, \quad \pno{p}{2} := \pno{q}{2} + \varepsilon_{2},
\]
where each $\varepsilon_{j}$ is a  realization of the Gaussian zero-mean i.i.d. random vector with standard deviation $10^{-4}$.
We average the results over $20$ realizations of the noise vector.

\begin{table}[ht!]
\parbox{.68\linewidth}{
\caption{Optimal approximations of the methods}
\label{tab:res_ex2b}
\centering
{\small
\begin{filecontents*}{res_ex2b_matlab.txt}
	k	n	d	LRA	VP$_{g}$(5)	VP$_{g}$	FASTGCD	FASTGCD(1)
   1.0000000e+00   5.4000000e+01   5.0000000e+01   2.1510822e-01   8.6818563e-02   8.3288822e-02   7.5080344e-01   8.9047151e+01
   2.0000000e+00   1.0400000e+02   1.0000000e+02   1.7086574e-01   3.1865444e-02   1.8096222e-02   6.5411626e-01   1.5310119e+01
   3.0000000e+00   2.0400000e+02   2.0000000e+02   1.0885681e-01   2.1766204e-02   2.0600790e-02   7.3549955e-01   1.3850185e+01
   4.0000000e+00   5.0400000e+02   5.0000000e+02   7.3436697e-02   3.7769180e-02   3.7711892e-02   2.7594067e-01   2.6427067e+00
   5.0000000e+00   1.0040000e+03   1.0000000e+03   8.3756351e-02   2.9060389e-02   2.4691850e-02   8.1224246e-01   3.3499732e+02
\end{filecontents*}
\pgfplotstabletypeset[
every head row/.style={after row=\hline},
columns={{d},{LRA},{VP$_{g}$},{VP$_{g}$(5)},FASTGCD},
every first column/.style={column type/.add={}{|}}]{res_ex2b_matlab.txt}
}}
\quad
\parbox{.27\linewidth}{
\caption{Avg. \# of iterations}
\label{tab:iters_ex2b}
\centering
{\small
\begin{filecontents*}{iters_ex2b_matlab.txt}
	k	n	d	LRA	VP$_{g}$(5)	VP$_{g}$	FASTGCD	FASTGCD(1)
   1.0000000e+00   5.4000000e+01   5.0000000e+01  -1.0000000e+00   4.9000000e+00   1.2350000e+01   4.5000000e+00   1.0000000e+00
   2.0000000e+00   1.0400000e+02   1.0000000e+02  -1.0000000e+00   4.9500000e+00   9.1000000e+00   5.8500000e+00   1.0000000e+00
   3.0000000e+00   2.0400000e+02   2.0000000e+02  -1.0000000e+00   4.9000000e+00   9.8000000e+00   5.8500000e+00   1.0000000e+00
   4.0000000e+00   5.0400000e+02   5.0000000e+02  -1.0000000e+00   4.7500000e+00   8.0500000e+00   5.7000000e+00   1.0000000e+00
   5.0000000e+00   1.0040000e+03   1.0000000e+03  -1.0000000e+00   4.8500000e+00   1.1000000e+01   4.3000000e+00   1.0000000e+00
\end{filecontents*}
\pgfplotstabletypeset[
every head row/.style={after row=\hline},
columns={{VP$_{g}$},FASTGCD},
every first column/.style={column type/.add={|}{}}]{iters_ex2b_matlab.txt}
}}
\end{table}

In Tables~\ref{tab:res_ex2b} and~\ref{tab:iters_ex2b}, we provide the approximation errors and numbers of iterations. For ``FASTGCD'' the average number of iterations is close to $5$. We also provide in Table~\ref{tab:res_ex2b} the results for  ``VP$_{g}$'' with number of iterations limited to $5$ (denoted by  ``VP$_{g}$(5)''). In this case, again ``VP$_{g}$''  achieves better approximation error for same number of iterations as ``FASTGCD''.  

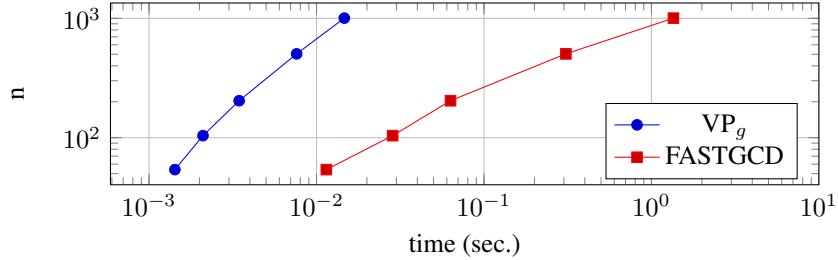
\begin{figure}[ht!]
\centering
\begin{tikzpicture}
\begin{loglogaxis}[ylabel=n,xlabel=time (sec.),xmax=10,legend pos = south east, height = 4cm, width = 11cm, grid=major]
\pgfplotstableread{times_ex2b_matlab.txt}\loadedtable %
\addplot table[y index=1,x index=3] {\loadedtable};
\addplot table[y index=1,x index=4] {\loadedtable};
\legend{VP$_{g}$,FASTGCD}
\end{loglogaxis}
\end{tikzpicture}
\caption{Comparison of the time per iteration, large $d$} 
\label{fig:speed2b}                                 
\end{figure}

It can be seen from Fig.~\ref{fig:speed2} and Fig.~\ref{fig:speed2b} , the time growth resembles $O(n)$ for the variable projection method, which confirms the results of Section~\ref{sec:ls_mosaic_relation}. The time growth for the iterations of ``FASTGCD'' resembles $O(n^2)$, which is consistent with complexity results of \cite{Bini.Boito10-fast}. Note that in Fig.~\ref{fig:speed2} the time needed for one iteration of the local optimization is of the same order as the total time reported in \cite{Bini.Boito10-fast} (including initial approximation), and therefore cannot be neglected.

\subsection{Example with several polynomials}\label{sec:ex_several_poly}
We consider the example of three polynomials \cite[Example 21.]{Christou.etal10ANM-ERES}.
\[
\begin{split}
\pno{p}{1}(s) = & -16.316 s^{11} + 182.73s^{10} -185.83 s^{9} +106.68 s^{8} - 266.22s^7 + 125.80 s^{6} \\
                & - 195.53 s^{5} + 243.81 s^4 + 23.013 s^3  + 64.186 s^2 -24.300 s - 43.810, \\
\pno{p}{2}(s)  = & 4.6618 s^{11} - 52.209 s^{10} + 53.094s^{9} - 30.481 s^{8} +76.064 s^{7} - 35.944s^{6} \\
                 &+ 55.866 s^5 - 69.659 s^4 - 6.5751s^3 -18.339s^2 +6.9428 s + 12.517, \\
\pno{p}{3}(s)  = & -4.1155s^{11} + 47.507 s^{10} - 59.034 s^9 + 2.2157 s^8 -45.276 s^7 +83.932 s^6 \\
                 &- 34.013 s^5 + 15.007 s^4  +4.3083 s^3 - 9.0031 s^2 + 14.297s -14.783.
\end{split}                 
\]
We are interested in the common divisor of degree $2$.

\begin{figure}
\centering
\includegraphics[height=6cm]{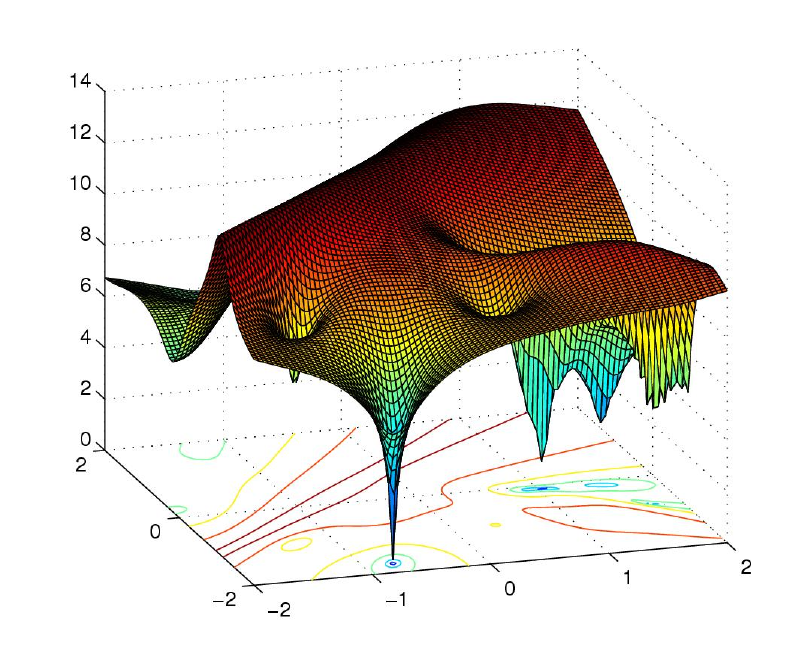}
\caption{The logarithm of the cost function: $\ln f^{(1)} ([-1 \; \widehat{a} \; \widehat{b}  ]^{\top})$.} \label{fig:cost_ex3}
\end{figure}

Since the degree of the common divisor is small, we will use the image representation\footnote{Optimization with variable projection in the Sylvester low-rank approximation is not applicable here, see the discussion in Section~\ref{sec:sing_ker}.}  and variable projection with respect to the common divisor (see Example~\ref{ex:varpro_h}), and solve the optimization problem \eqref{eq:outer} (minimize the function $f_1(\widehat{h})$ over $\pspace{d} \setminus \{0\}$). Since $f_1(\widehat{h})$ is invariant of scaling of the parameter, this is a problem of optimization on a projective space. We fix a coordinate chart in this space and optimize only over the polynomials $\widehat{h}(z) = \widehat{b} z^2 + \widehat{a} z -1$. 
In Fig.~\ref{fig:cost_ex3}, we plot the cost function $f^{(1)}(\bmx -1 & \widehat{a} & \widehat{b}\emx^{\top})$ evaluated on a $100 \times 100$ grid in the box $[-2;2] \times [-2;2]$. In Fig.~\ref{fig:cost_ex3}, we see that $f^{(1)}$ possesses many local minima and a large Lipshitz constant.

We consider several initial approximations. The polynomials  $\widehat{h}_{kl0}$, where, $1 \le k < l \le 3$, are the initial approximations obtained from Sylvester sub-resultant matrices of two polynomials $\pno{p}{k}$ and $\pno{p}{l}$ by Algorithm~\ref{alg:iniapprox_sylv_subr}. The polynomial $\widehat{h}_{1230}$ denotes the result of Algorithm~\ref{alg:iniapprox_sylv_subr} for the structure $\sylv{d}(\bfp)$, and $\widehat{h}_{l23f}$ denotes the result of the Algorithm~\ref{alg:iniapprox_sylv_subr} for the structure $\sylv{d}^{(full)}(\bfp)$. $\widehat{h}_{123m}$ denotes the result of Algorithm~\ref{alg:iniapprox_sylv_ext} (for $K = 32$). Finally, $\widehat{h}_{ref}$ denotes the result provided in \cite[Example 21.]{Christou.etal10ANM-ERES}:
\[
c(z^2  -11.28371806974011z + 11.64469379842480),
\]
where the constant $c$ is chosen to conform to the  normalization.
All the polynomials are normalized to be of the form $\widehat{h}(z) = \widehat{b} z^2 + \widehat{a} z -1$.

From all initial approximations, we run the optimization algorithm for \eqref{eq:outer} (with the maximum number of iterations $300$). The initial approximations ($\widehat{a}_0, \widehat{b}_0$), the initial value of the cost function ($f_0$), the distance to the reference polynomial ($\rho_{ref}$), the point of the local minimum ($\widehat{a}_{opt}, \widehat{b}_{opt}$), the cost function value at the minimum ($f_{opt}$) and the number of iterations needed (iter.) are shown in Table~\ref{tab:tab3res}.

\begin{table}
\caption{Optimal approximations of the methods}
\label{tab:tab3res}
\centering
{\small 
\begin{filecontents*}{res_ex3.txt}
   9.3343837e-01  -8.0410519e-01   1.5242294e+00   9.3279472e-01   9.3596812e-01   9.6912949e-01   9.6900084e-01
  -7.1535075e-02  -1.8681338e+00  -2.4949472e-01   1.6918860e-02   1.1406116e-02  -8.6458489e-02  -8.5876024e-02
   9.4732594e+00   6.9971253e+01   1.7733503e+02   5.4183288e+02   4.8617449e+02   1.3131040e-02   4.6901069e-04
   3.8345173e-02   2.5140302e+00   5.7883485e-01   1.0898473e-01   1.0273741e-01   5.9650486e-04   0.0000000e+00
   9.6919817e-01  -7.4130790e-01   1.5231236e+00   9.8649457e-01   9.8649451e-01   9.6919817e-01   9.6919817e-01
  -8.5982870e-02  -1.7419861e+00  -1.4068562e-01  -8.7686361e-02  -8.7686355e-02  -8.5982870e-02  -8.5982871e-02
   3.6354379e-07   1.4900084e-06   1.1253553e+01   8.2990802e-04   8.2990802e-04   3.6354379e-07   3.6354379e-07
   5.0000000e+00   4.0000000e+00   4.0000000e+00   7.0000000e+00   7.0000000e+00   2.0000000e+00   2.0000000e+00
\end{filecontents*}
\pgfplotstabletypeset[
	create on use/value/.style={ create col/set list={$\widehat{a}_0$,$\widehat{b}_0$,$f_0$, $\rho_{ref}$, $\widehat{a}_{opt}$,$\widehat{b}_{opt}$,$f_{opt}$, iter.}}, 
	columns/value/.style={string type},
every head row/.style={after row=\hline},
columns/0/.style={column name={$\widehat{h}_{130}$}, column type/.add={@{\,}|@{\,}}{}},
columns/1/.style={column name={$\widehat{h}_{230}$}, column type/.add={@{\,}|@{\,}}{}},
columns/2/.style={column name={$\widehat{h}_{120}$}, column type/.add={@{\,}|@{\,}}{}},
columns/3/.style={column name={$\widehat{h}_{1230}$}, column type/.add={@{\,}|@{\,}}{}},
columns/4/.style={column name={$\widehat{h}_{123f}$}, column type/.add={@{\,}|@{\,}}{}},
columns/5/.style={column name={$\widehat{h}_{123m}$}, column type/.add={@{\,}|@{\,}}{}},
columns/6/.style={column name={$\widehat{h}_{ref}$}, column type/.add={@{\,}|@{\,}}{}},
columns ={value,0,1,2,3,4,5,6},
]{res_ex3.txt}}
\end{table}

In Table~\ref{tab:tab3res}, we see that $\widehat{h}_{130}$ and $\widehat{h}_{123m}$ give the best answer. The polynomial are $\widehat{h}_{123m}$ is the closest to the reference polynomial  \cite[Example 21.]{Christou.etal10ANM-ERES}, and gives the same result as optimization started from the reference polynomial.
We also see from Table~\ref{tab:tab3res}, that the initial approximation $\widehat{h}_{123f}$ obtained from  improved the Sylvester subresultant $\sylv{d}^{(full)}(\bfp)$ (suggested in \cite{Agrawal.etal04SP-Common}) is slightly better (and closer to the reference polynomial) than $\widehat{h}_{1230}$ (obtained from the subresultant $\sylv{d}(\bfp)$). Both give a good solution, but not the optimal one (they fall into a neighboring local minimum).

We also see that the methods converged to different local minima,  shown in Fig.~\ref{fig:cost_ex3}. Note that Fig.~\ref{fig:cost_ex3} does not reflect the values of the cost function at local minima, for example, the global minimum (for $\widehat{h}_{130}$ and $\widehat{h}_{123m}$) is not visible in Fig.~\ref{fig:cost_ex3}. This is explained by the fact that Fig.~\ref{fig:cost_ex3} is evaluated on a grid and some minima may not be captured by the grid. This also shows  intrinsic complexity of the optimization problem.

\subsection{Applicability of the kernel representation for \texorpdfstring{$N>2$}{N > 2}: singularity of \texorpdfstring{$\Gamma$}{Gamma} matrix}\label{sec:sing_ker}
In this subsection we show that on a particular example of $N = 3$ polynomials, in the Sylvester low-rank approximation approach the corresponding $\Gamma$ matrix from \cite{Usevich.Markovsky13JCAM-Variable} is essentially singular. Consider three polynomials $\pno{p}{1},\pno{p}{2},\pno{p}{3} \in \pspace{2}$, and $d = 1$. Then the corresponding generalized Sylvester subresultant matrix is
\[
\sylv{1}^{(1)} (\bfp) =
\bmx
\multmat{\pno{p}{2}}{1}  & -\multmat{\pno{p}{1}}{1}   &          \\
\multmat{\pno{p}{3}}{1}  &                               &  -\multmat{\pno{p}{1}}{1} \\
\emx.
\]
By Proposition~\ref{prop:sylv_mosaic} and the results of \cite{Usevich.Markovsky13JCAM-Variable}, we have that the corresponding $\Gamma$ matrix  has the form  $\Gamma(\widehat{\bfu}) = G(\bfu) G^{\top} (\widehat{\bfu})$, where 
\[
G(\bfu) = 
\bmx
-\multmat{\pno{\widehat{u}}{2}}{2}  & \multmat{\pno{\widehat{u}}{1}}{2}   &          \\
-\multmat{\pno{\widehat{u}}{3}}{2}  &                               &  \multmat{\pno{\widehat{u}}{1}}{2} \\
\emx,
\]
and $\widehat{\bfu} \in \pspace{[1 1 1]}$ such that  $\sylv{1}^{(1)} \widehat{\bfu} = 0$.
It can be easily checked that the polynomial matrix $\Gamma(\widehat{\bfu})$ has (symbolic) determinant $0$. This is also confirmed by running the optimization method with the help of SLRA package.

\section{Conclusions}
We have developed methods based on the variable projection principle, for optimization in the direct parameterization and Sylvester low-rank approximation. The advantages of the developed methods are that they have proven complexity results, and have available implementation that allows to use different optimization methods and different stopping criteria. 

The methods for optimization in the direct  parameterization have linear complexity in the degrees of the polynomials if the degree of the common divisor $d$ is small or if $d$ is large. The methods provide accurate results matching the accuracy of other existing methods. We also showed that the methods based on direct  parameterization perform better than the methods based on the kernel representation. The latter have higher computational complexity and have issues of intrinsic singularity of the $\Gamma$ matrix for $N > 2$ and ill-conditioning of $\Gamma$ when the degree $d$ is small. 


\section*{Acknowledgement}
The research leading to these results has received funding from the European Research Council under the European Union's Seventh Framework Programme (FP7/2007-2013) / ERC Grant agreement no. 258581 ``Structured low-rank approximation: Theory, algorithms, and applications'' and Grant Agreement No. 320594 DECODA project.
We also thank the anonymous reviewers for valuable suggestions, which led to significant improvements in the presentation of the paper.

\appendix
\section{Basic properties of the ACD problem}\label{sec:acd_basic}

\begin{lem}\label{lem:mod_class_relation}
For the sets $\pspacegcd{d}$ defined in Sections~\ref{sec:prob_stat} and~\ref{sec:prob_stat}
\[
\pspacegcd{d} =
\begin{cases}
\pspacecd{d}, & \mbox{if }  \bbF = \bbC \mbox{ or } (\bbF = \bbR \mbox{ and } d \mbox{ is even}),\\
\pspacecd{d} \cup \pspacecd{d+1}, & \mbox{if }  \bbF = \bbR \mbox{ and } d \mbox{ is odd}.
\end{cases}
\]
\end{lem}
\begin{pf}
Evidently, $\pspacecd{d} \subset \pspacegcd{d}$, $\pspacecd{d+1} \subset \pspacegcd{d}$ and $0 \in \pspacecd{d}$. Let $\bfp= (\pno{p}{1},\ldots,\pno{p}{N}) \in \pspacegcd{d} \setminus \{ 0 \}$ and $h \in \gcd (\bfp)$, where $\deg h = d' > d$. Then, $h$ can be factorized as 
\[
h(z) = \sum\limits_{k=1}^{d'} (z - \alpha_k),
\]
where $\alpha_k \in \bbC \cup \{ \infty \}$. If $\bbF = \bbC$ then the $\pno{p}{k}$ have a common divisor of degree $d$ and $\bfp \in \pspacecd{d}$. If $\bbF = \bbR$, then every $\alpha_k$ has its conjugate counterpart. Hence, $\bfp$ have common divisors of degrees $2l$ for any $l: 0 \le 2l \le d'$. Therefore, $\bfp \in \pspacecd{d}$ if $d$ is even and $\bfp \in \pspacecd{d+1}$ if $d$ is odd.
\end{pf}

\begin{lem}\label{lem:mod_class_closed}
For any $0 \le d \le n_{min}$,  sets $\pspacecd{d}, \pspacegcd{d}$ are closed subsets of $\pspacen$. 
\end{lem}
\begin{pf}
Denote by $\calS_d \subset \pspace{d}$ the polynomials with $2$-norm equal to $1$.
Then $\pspacecd{d}$ is the image of the infinitely smooth map 
$\pspace{\bfn-d} \times \calS_d \to \pspace{\bfn}$ defined as 
$(\widehat{\bfg}, \widehat{h})\mapsto \widehat{\bfg} \cdot \widehat{h}$.
Since the domain of the definition is closed and the map is continuous, $\pspacecd{d}$ is closed.
Finally by Lemma~\ref{lem:mod_class_relation},  $\pspacegcd{d}$ can be expressed as a union of at most two sets  $\pspacecd{d_1}$ and $\pspacecd{d_2}$, and therefore it is also closed.
\end{pf}

\section{Properties of the extended Sylvester sub resultants}\label{sec:sylv_full_properties}
The following lemma is well-known in the computer algebra community \cite{Rupprecht99JPAA-algorithm,Kaltofen.etal06conf-Approximate}. We present it in a modified version, which takes care of possible zero polynomials.
\begin{lem}[{\cite[Lemma 2.1, adjusted]{Kaltofen.etal06conf-Approximate}}]\label{lem:fracs}
Let $\bfp \in \pspacen$, $\pno{u}{1} \in \pspace{n_1-d} \setminus \{0\}$ and $\pno{u}{2} \in \pspace{n_2-d}$, $\ldots$, $\pno{u}{N} \in \pspace{n_N-d}$ such that
\begin{equation}\label{u_k}
\pno{u}{k} \pno{p}{1} - \pno{u}{1} \pno{p}{k} = 0, \;\mbox{where} \; \pno{u}{k} \in \pspace{n_k},
\; \forall k = 2, \ldots, N.
\end{equation}
Then the polynomials $\pno{p}{1}, \ldots, \pno{p}{N}$ have a common divisor of degree at least $d$.
\end{lem}
\begin{pf}
If $\pno{p}{1} = 0$, then all polynomials $\pno{p}{k}$ are zero by \eqref{u_k}.
Then we are left to consider $\pno{p}{1} \neq 0$. (The case when all other polynomials are zero is trivial. )Without loss of generality assume that there exists $2 \le K \le N$ such that $\pno{p}{k} \neq 0$ for all $2 \le k \le N$. Then we need to prove that $\deg \gcd (\pno{p}{1}, \ldots, \pno{p}{K}) \ge d$.

Since $\pno{u}{1} \pno{p}{k} \neq 0$, we have that $\pno{u}{k} \neq 0$ for $k \ge 2$. Let us rewrite \eqref{u_k} as
\begin{equation}\label{eq:fracs}
\frac{\pno{p}{1}}{\pno{u}{1}} = \frac{\pno{p}{2}}{\pno{u}{2}} = \cdots = \frac{\pno{p}{K}}{\pno{u}{K}} = \frac{a}{b},
\end{equation}
where $a/b$ is an irreducible fraction. Then we have that
\begin{equation}
\pno{p}{k} = \frac{\pno{u}{k} a}{b},
\label{eq:pk_uk}
\end{equation}
and $\pno{u}{k}$ should be divisible by $b$. Therefore, $\pno{p}{k}$ have a common divisor $a$, with $\deg a \ge d$, which can be established by counting dimensions in \eqref{eq:fracs}.
\end{pf}

From Lemma~\ref{lem:fracs}, Lemma~\ref{lem:standard_extended_sylvester} easily follows. 

\begin{pf}[{Proof of Lemma~\ref{lem:standard_extended_sylvester}}]
$\boxed{\Leftarrow}$ Suppose that $\sylv{d}$ is rank deficient. Then there exists 
$\bfu = (\pno{u}{1},\ldots,\pno{u}{N}) \in \pspace{\bfn-d} \setminus \{\bfzero\}$ such that $\sylv{d} ({\bfp}) \bfu = \bfzero$. Therefore, the equations \eqref{u_k} are satisfied. Since $\pno{p}{1} \neq 0$ and $\bfu \neq 0$, we have that $\pno{u}{1}$. Therefore, by Lemma~\ref{lem:fracs}, we have that $\bfp \in \pspacegcd{d}$.

$\boxed{\Rightarrow}$ Since $\bfp \in \pspacegcd{d}$, there exist polynomials 
$\pno{u}{k} \in \pspace{n_k - d}$ and $h \in \pspace{d}$ such that $\pno{p}{k} = \pno{u}{k} h$.
Then, immediately,  the equations \eqref{u_k} are satisfied. Therefore, for the vector 
$\bfu = (\pno{u}{1},\ldots,\pno{u}{N}) \in \pspace{\bfn-d} \neq 0$, we have that
$\sylv{d} ({\bfp}) \bfu = \bfzero$. 
\end{pf}

Now we are in a position to prove Proposition~\ref{prop:sylv_equiv}.

\begin{pf}[Proof of Proposition~\ref{prop:sylv_equiv}]
The ``only if'' part is trivial. Indeed, one can construct $\bfu$, as in Lemma~\ref{lem:fracs}.
In order to prove the ``if'', denote by $\sylv{d}^{(1)} = \sylv{d}$, and denote by $\sylv{d}^{(k)}$, $k \ge 2$, shifted Sylvester subresultants
\[
\sylv{d}^{(k)}(\pno{p}{1}, ,\ldots, \pno{p}{N}) \defeq
\sylv{d}(\pno{p}{k}, \pno{p}{1}, \ldots, \pno{p}{k-1}, \pno{p}{k+1},\ldots, \pno{p}{N})
\bmx
  &   & I_{l_k} & \\
  & I_{m_1} &   & \\
I_{l_1} &   &   & \\
  &   &   & I_{m_2}
\emx,
\]
where $l_k \defeq n_k -d+1$, $m_1 := \elsum{(l_{2} , \ldots,l_{k-1})}$ and $m_2:= \elsum{(l_{k+1}, \ldots,l_{N})}$. For example, for three polynomials we have that
\[
\begin{split}
&\sylv{d}^{(2)}(\pno{p}{1}, \pno{p}{2}, \pno{p}{3}) \defeq
\bmx
-\multmat{\pno{p}{2}}{n_1-d}  & \multmat{\pno{p}{1}}{n_2-d}   &          \\
                              & \multmat{\pno{p}{3}}{n_2-d}  &  -\multmat{\pno{p}{2}}{n_3-d} \\
\emx, \\
&\sylv{d}^{(3)}(\pno{p}{1}, \pno{p}{2}, \pno{p}{3}) \defeq
\bmx
 -\multmat{\pno{p}{3}}{n_1-d} &                              & \multmat{\pno{p}{1}}{n_3-d}  \\
                              & -\multmat{\pno{p}{3}}{n_2-d} & \multmat{\pno{p}{2}}{n_3-d} \\
\emx. 
\end{split}
\]
Note that any matrix $\sylv{d}^{(k)} (\bfp)$ can be extracted from the matrix $\sylv{d}^{(full)} (\bfp)$ by selecting corresponding block rows (and, possibly, negation). Therefore, if $\bfu \in \pspace{\bfn -d}$ is a right annihilating vector $\bfu \in \pspace{\bfn -d}$ of $\sylv{d}^{(full)} (\bfp)$, it is also annihilating vector of $\sylv{d}^{(k)} (\bfp)$. Let us select $k$ such that 
$\pno{u}{k} \neq 0$. Then we have that
\[
\sylv{d} (\pno{p}{k}, \pno{p}{1}, \ldots, \pno{p}{k-1}, \pno{p}{k+1},\ldots, \pno{p}{N})
\vcol(\pno{u}{k}, \pno{u}{1}, \ldots, \pno{u}{k-1}, \pno{u}{k+1},\ldots, \pno{u}{N}) = 0,
\]
and by Lemma~\ref{lem:fracs}, the polynomials have $\deg\gcd \ge d$.
\end{pf}

\section{Least-squares and least-norm problems}\label{sec:ls_ln}
In this paper, we use the duality between least-squares and least-norm problems. Next, we give an overview of these problems.

\begin{prob}[Weighted least-squares problem]
Let $A \in \bbR^{m \times n}$, $m \ge n$, $b \in \bbR^{m}$, and $v \in [0,\infty)^{m}$, such that $\rrank \diag(\sqrt{v})  A = n$
\begin{equation}
\minim_{x \in \bbR^{n}} \|A x - b\|^2_v.
\label{eq:lls}
\end{equation}
\end{prob}
The problem  \eqref{eq:lls} is the orthogonal projection of $b$ on the image of $A$ in the seminorm $\|\cdot\|_v$. The solution can be found by rewriting the cost function in \eqref{eq:lls}  as  $\|\lsfun{LS}\|^2_2$, where $\lsfun{LS} = \diag(\sqrt{v})(b - A x)$. Then the solution of \eqref{eq:lls} is 
\begin{equation}
\begin{split}
&x_* = (A^{\top} \diag(v) A)^{-1} A^{\top} \diag(v) b, \\ 
&\lsfunopt{LS} = \diag(\sqrt{v}) b -  diag(\sqrt{v}) A (A^{\top} \diag(v) A)^{-1} A^{\top} \diag(v)  b,\\
&\optfun{LS} = \|\lsfunopt{LS}\|_2^2 = \|b\|^2_v -   b^{\top} \diag(v)  A (A^{\top} \diag(v) A)^{-1} A^{\top} \diag(v) b.\\
\end{split}\label{eq:lls_sol}
\end{equation}
The least-squares problem  \eqref{eq:lls} is closely connected to the following dual problem:
\begin{prob}[Weighted least-norm problem]
Let $A \in \bbR^{m \times n}$, $m \ge n$, $c \in \bbR^{m}$, and $w \in (0,\infty]^{m}$, such that $\rrank \diag(\sqrt{w^{-1}}) A = n$
\begin{equation}
\minim_{z \in \bbR^{m}} \|c - z \|^2_{w} \quad\sto\quad A^{\top} z = 0.
\label{eq:lln}
\end{equation}
\end{prob}
The problem \eqref{eq:lln} is to find the orthogonal projection of the vector $c$ on the kernel of the matrix $A$ in the norm $\|\cdot\|_{w}$. Changing variables as $c - z = \diag(\sqrt{w^{-1}}) \lsfun{LN}$ The cost function can be rewritten as $\|\lsfun{LN}\|^2_2$, and the constraint $A^{\top}z = 0$ as 
\[
A^{\top} \diag(\sqrt{w^{-1}}) \lsfun{LN}  = A^{\top} c.
\]
Therefore, the solution of \eqref{eq:lln} is the following
\begin{equation}
\begin{split}
&z_* = c - \diag(w^{-1}) A (A^{\top} \diag(w^{-1}) A)^{-1} A^{\top} c, \\
&\lsfunopt{LN} = \diag(\sqrt{w^{-1}}) A (A^{\top} \diag(w^{-1}) A)^{-1} A^{\top} c, \\
&\optfun{LN} = \|\lsfunopt{LN}\|_2^2  = c^{\top} A (A^{\top} \diag(w^{-1}) A)^{-1} A^{\top} c.
\end{split}\label{eq:lln_sol}
\end{equation}
One can see that the expressions for the solutions of \eqref{eq:lls} and \eqref{eq:lln} have a similar form. In particular, if we have $c = \diag(w^{-1})b$ and $v = w^{-1}$, then
\begin{equation}
\begin{split}
&\optfun{LS} = \|b\|^2_v - \optfun{LN}, \\
&\lsfunopt{LS} = \diag(\sqrt{v}) b - \lsfunopt{LN}.
\end{split}
\label{eq:lls_lln_sol}
\end{equation}

\bibliographystyle{unsrt}
\section*{References}
{\small
\bibliography{agcd}
}

\end{document}